\documentclass[titlepage,final,11pt]{article}
\usepackage{graphics}
\usepackage{graphicx}
\usepackage{subfigure}
\usepackage{epsfig}
\usepackage{enumerate}
\usepackage{amsfonts}
\usepackage{theorem}
\usepackage{array}
\usepackage[reqno]{amsmath}
\usepackage{color}

\theoremstyle{change}
\newtheorem{proclaim}{PROCLAIM}[section]
\newtheorem{theorem}[proclaim]{Theorem}

\newtheorem{lemma}[proclaim]{Lemma}
\newtheorem{proposition}[proclaim]{Proposition}

\def\state #1. { \noindent{\bf#1.\enspace}}

\newcommand{\comp}{\,{\raise 1pt \hbox{$\scriptstyle\circ$}}\,}

\newcommand{\grill}{{\scriptscriptstyle\#}}

\newcommand{\dist}{\mathop{\rm dist}\nolimits}
\newcommand{\exs}{\mathop{\rm exs}\nolimits}

\newcommand{\hypo}{\mathop{\rm hypo}}

\newcommand{\cl}{\mathop{\rm cl}}

\newcommand{\nt}{\mathop{\rm int}}
\newcommand{\lev}{\mathop{\rm lev}}

\newcommand{\uscfcns}{\mathop{\textrm{usc-fcns}}}

\newcommand{\reals}{{I\kern-.35em R}}
\newcommand{\Reals}{\overline{\reals}}

\newcommand{\natnums}{{{\rm l} \kern -.13em {\rm N} }}
\newcommand{\nats}{{I\kern -.35em N}}
\newcommand{\snats}{{I\kern -.29em N}}
\newcommand{\rats}{{Q\kern -.64em \raise 1pt \hbox{$\scriptstyle |$}\;\,}}
\newcommand{\srats}
    {{Q\kern -.56em \raise 1.2pt \hbox{$\scriptscriptstyle /$}\,}}
\newcommand{\ints}{Z\kern -.46em Z}
\newcommand{\ball}{{I\kern -.35em B}}
\newcommand{\Ex}{\mathbb{E}}
\newcommand{\pluss}{\hskip1pt \raise1pt\vbox{\hrule width6pt \vskip1pt \hrule
                    width6pt} \kern-4pt{\lower1pt\hbox{\vrule height6pt
            \kern1pt\vrule height6pt}}\hskip5pt}
\newcommand{\eop}
    {\hfill{$\vcenter{\hrule height1pt \hbox{\vrule width1pt height5pt
     \kern5pt \vrule width1pt} \hrule height1pt}$} \medskip}

\newcommand{\setd}{{ d \kern -.15em l}}
\newcommand{\hatsetd}{ d \hat{\kern -.15em l }}

\renewcommand{\epsilon}{\varepsilon}
\renewcommand{\phi}{\varphi}

\newcommand{\tto}{\;{\lower 1pt \hbox{$\rightarrow$}}\kern -12pt
           \hbox{\raise 2.5pt \hbox{$\rightarrow$}}\;}
\newcommand{\overto}[1]{\,{\raise 0pt\hbox{$\rightarrow$}}\kern -9pt
     \hbox{\lower 3pt \hbox{$\scriptscriptstyle#1$}}\hskip6pt}
\newcommand{\underto}[1]{\,{\lower 1pt\hbox{$\rightarrow$}}\kern -9pt
     \hbox{\raise 4pt \hbox{$\,\scriptscriptstyle#1$}}\hskip7pt}
\newcommand{\bigoverto}[1]{{\raise 0pt\hbox{$\,\longrightarrow$}}\kern -16pt
     \hbox{\lower 3pt \hbox{$\scriptscriptstyle#1$}}\hskip4pt}
\newcommand{\bigunderto}[1]{\,{\lower 1pt\hbox{$\longrightarrow$}}\kern -16pt
     \hbox{\raise 4pt \hbox{$\,\scriptscriptstyle#1$}}\hskip6pt}
\newcommand{\bigbigto}[2]{\,{\raise 0pt\hbox{$\,\longrightarrow$}}\kern -16pt
     \hbox{\lower 3pt \hbox{$\scriptscriptstyle#2$}}\kern -10pt
     \hbox{\raise 4pt \hbox{$\,\scriptscriptstyle#1$}}\hskip7pt}
\newcommand{\downto}{{\raise 1pt \hbox{$\scriptstyle \,\searrow\,$}}}
\newcommand{\upto}{{\raise 1pt \hbox{$\scriptstyle \,\nearrow\,$}}}

\newcommand{\notimply}
    {\quad\hbox{$\Longrightarrow \kern -14pt {/}$}\hskip6pt\quad}

\newcommand{\nLim}{\mathop{\rm Lim}\nolimits}

\newcommand{\nInnLim}{\mathop{\rm InnLim}\nolimits}
\newcommand{\nOutLim}{\mathop{\rm OutLim}\nolimits}

\newcommand{\low}[1]{{\lower1pt \hbox{$\scriptstyle #1$}}}
\newcommand{\loww}[1]{{\lower2pt \hbox{$\scriptstyle #1$}}}
\newcommand{\high}[1]{{\raise1pt \hbox{$\scriptstyle #1$}}}

\newcommand{\cA}{{\cal A}}
\newcommand{\cB}{{\cal B}}

\newcommand{\cN}{{\cal N}}

\newcommand{\cR}{{\cal R}}

\newcommand{\nlev}{\mathop{\lev}\nolimits}

\newcommand{\nlim}{\mathop{\rm lim}\nolimits}
\newcommand{\nliminf}{\mathop{\rm liminf}\nolimits}
\newcommand{\nlimsup}{\mathop{\rm limsup}\nolimits}
\newcommand{\ninf}{\mathop{\rm inf}\nolimits}
\newcommand{\nsup}{\mathop{\rm sup}\nolimits}
\newcommand{\nmax}{\mathop{\rm max}\nolimits}
\newcommand{\nmin}{\mathop{\rm min}\nolimits}

\newcommand{\nargmin}{\mathop{\rm argmin}\nolimits}

\newcommand{\bfxi}{\mbox{\boldmath $\xi$}}

\newcommand{\bfx}{\mbox{\boldmath $x$}}
\newcommand{\bfy}{\mbox{\boldmath $y$}}

\newcommand{\bfz}{\mbox{\boldmath $z$}}

\newcommand{\lwdy}[2]{\mathrel{\mathop
        {\raisebox{0.1ex}{\null$#1$}}{\hbox{\kern -1.0em
    {\raisebox{-0.8ex}{$\scriptstyle{\;\to #2}$}}}}}}
\newcommand{\lwwdy}[2]{\mathrel{\mathop
        {\raisebox{0.2ex}{\null$#1$}}{\hbox{\kern -1.0em
    {\raisebox{-1.1ex}{$\scriptstyle{\;\to #2}$}}}}}}
\newcommand{\slwwdy}[2]{\scriptsize{{\mathrel{\mathop
        {\raisebox{0.2ex}{\null$#1$}}{\hbox{\kern -1.0em
    {\raisebox{-1.1ex}{$\scriptstyle{\;\to #2}$}}}}}}}}

\def\eto{\,{\lower 1pt\hbox{$\rightarrow$}}\kern -11pt
     \hbox{\raise 4pt \hbox{$\, \scriptstyle e$}}\hskip7pt}
\def\hto{\,{\lower 1pt\hbox{$\rightarrow$}}\kern -11pt
     \hbox{\raise 4pt \hbox{$\, \scriptstyle h$}}\hskip7pt}
\def\pto{\,{\lower 1pt\hbox{$\rightarrow$}}\kern -11pt
     \hbox{\raise 4.5pt \hbox{$\, \scriptstyle p$}}\hskip7pt}
\def\cto{\,{\lower 1pt\hbox{$\rightarrow$}}\kern -11pt
     \hbox{\raise 4pt \hbox{$\, \scriptstyle c$}}\hskip7pt}
\def\Nto{\,{\lower 1pt\hbox{$\rightarrow$}}\kern -12pt
     \hbox{\raise 4pt \hbox{$\, \scriptstyle N$}}\hskip7pt}

\def\dy#1\until#2{\mathrel{\mathop
        {\null#1}\limits^{\hbox{\lower 1.3ex \hbox{$\scriptstyle{\;\to}$}}}}
        {\hbox{{\kern -.2em${\null}
       ^{\hbox{\raise .2ex \hbox{$\scriptstyle{#2}$}}}$}}}}
\def\hidy#1\until#2{\mathrel{\mathop
        {\null#1}\limits^{\hbox{\lower 1.0ex \hbox{$\scriptstyle{\;\to}$}}}}
        {\hbox{{\kern -.2em${\null}
	^{\hbox{\raise 0.9ex \hbox{$\scriptstyle{#2}$}}}$}}}\!}
\def\lody#1\until#2{\mathrel{\mathop
        {\null#1}\limits_{\hbox{\raise 0.5ex \hbox{$\scriptstyle{\to}$}}}}
        {\hbox{{\kern -.25em${\null}
       _{\hbox{\lower 0.7ex \hbox{$\scriptstyle{#2}$}}}$}}}}
\def\lowdy#1\until#2{\mathrel{\mathop
        {\null#1}\limits_{\hbox{\raise 1.0ex \hbox{$\scriptstyle{\to}$}}}}
        {\hbox{{\kern -.25em${\null}
       _{\hbox{\lower 0.9ex \hbox{$\scriptstyle{#2}$}}}$}}}}

\newcommand{\espl}{\textup{\textrm{e-spl}}}
\newcommand{\pafcns}{\textup{\textrm{pa-fcns}}}

\begin{document}


\begin{center}
\begin{large}
{\bf Approximations of Semicontinuous Functions\\ with Applications to
Stochastic Optimization and Statistical Estimation}
\smallskip
\end{large}
\vglue 0.7truecm
\begin{tabular}{c}
  \begin{large} {\sl Johannes O. Royset 
                                  } \end{large} \\
  Operations Research Department\\
  Naval Postgraduate School\\
  joroyset@nps.edu\\
\end{tabular}

\vskip 0.2truecm

\end{center}

\vskip 1.3truecm

\noindent {\bf Abstract}. \quad Upper semicontinuous (usc) functions arise in
the analysis of maximization problems, distributionally robust optimization,
and function identification, which includes many problems of nonparametric
statistics. We establish that every usc function is the limit of a
hypo-converging sequence of piecewise affine functions of the
difference-of-max type and illustrate resulting algorithmic possibilities in
the context of approximate solution of infinite-dimensional optimization
problems. In an effort to quantify the ease with which classes of usc
functions can be approximated by finite collections, we provide upper and
lower bounds on covering numbers for bounded sets of usc functions under the
Attouch-Wets distance. The result is applied in the context of stochastic
optimization problems defined over spaces of usc functions. We establish
confidence regions for optimal solutions based on sample average
approximations and examine the accompanying rates of convergence. Examples
from nonparametric statistics illustrate the
results.\\

\vskip 0.5truecm

\halign{&\vtop{\parindent=0pt
   \hangindent2.5em\strut#\strut}\cr
{\bf Keywords}: \ hypo-convergence, Attouch-Wets distance, approximation theory, solution stability, 
                         \hfill\break
\hglue 1.35cm stochastic optimization, epi-splines, rate of convergence.
                         \cr\cr

{\bf Date}:\quad \ \today \cr}

\baselineskip=15pt

\section{Introduction}

Extended real-valued upper-semicontinuous (usc) functions on $\reals^n$ are
fundamental in the study of finite-dimensional constrained maximization
problems as essentially all such problems can be represented by usc
functions. They arise in probability theory with distribution and
c\`{a}dl\`{a}g functions also being usc. Emerging applications of usc
functions in infinite-dimensional problems include nonparametric statistical
$M$-estimation \cite{RoysetWets.18}, distributionally robust optimization
\cite{RoysetWets.16b}, and more generally function identification
\cite{RoysetWets.15b}. In these applications, optimization problems are
formulated over spaces of usc functions. Regardless of the setting, it
becomes important to have means to approximate usc functions as well as an
understanding of the difficulty with such an undertaking. This article
provides three main results in these directions: (i) We establish that every
usc function is the limit of a hypo-converging sequence of mesh-free
piecewise affine functions of the difference-of-max type. Thus, as a
corollary, the difference-of-convex (dc) functions are hypo-dense in spaces
of usc functions. With the advances in computational treatment of dc
functions (see for example \cite{CuiPangSen.18}), this leads to numerous
algorithmic possibilities, which we illustrate in the context of function
identification problems. (ii) We provide upper and lower bounds on covering
numbers for bounded sets of usc functions under the Attouch-Wets (aw)
distance and thereby quantify the ease with which classes of usc functions
can be approximated by finite collections. (iii) For stochastic optimization
problems defined over spaces of usc functions, we establish confidence
regions for optimal solutions in terms of the aw-distance and sample average
approximations, with rates of convergence as the sample size grows. The
result requires only semicontinuity of the objective function and therefore
applies in challenging settings such as simulation optimization of
``black-box'' stochastic systems where little structure may be known.

The consideration of approximations in the sense of hypo-convergence, which
is metrized by the aw-distance, is natural and convenient in many
applications. If an usc function is approximated in this sense, then the
maximizers of the approximating function will be ``near'' those of the actual
function. This is exactly the desired property when the usc function
represents a constrained maximization problem. It is also the goal when the
usc functions is a probability density function and we need to estimate its
modes; the approximating density will have modes ``near'' the actual modes.
The situation is similar when the usc function is a surrogate model in an
engineering design problem that needs to be maximized to find an optimal
design; see Section 5 for an example. The notion of approximation is further
motivated in the context of distribution functions by the fact that for such
functions hypo-convergence is equivalent to convergence in distribution, a
property that is leveraged to address optimization under stochastic ambiguity
in \cite{RoysetWets.16b}. An alternative focus on approximations in the sense
of uniform convergence would have limited the scope to finite-valued
continuous functions with common compact domains, which is too restrictive in
many applications. Hypo-convergence permits treatment of usc functions
defined on any subset of $\reals^n$.

The study of hypo-converging usc functions and, in parallel, epi-converging
extended real-valued lower-semicontinuous (lsc) functions has a long history,
with important accomplishments in convex and nonsmooth analysis as well as
the approximation theory of maximization and minimization problems; see
\cite{VaAn} for details. Connections to probability theory are established in
\cite{SalinettiWets.86,SalinettiWets.86a} and more recently in
\cite{RoysetWets.16b}. The first formulation of infinite-dimensional
optimization problems over spaces of semicontinuous functions appears in
\cite{RoysetWets.14}, with theoretical developments in \cite{RoysetWets.15b}.
In particular, the latter reference defines the class of epi-splines (see
also \cite{Royset.18}), which are piecewise polynomial functions, and
establishes that the class is dense in spaces of semicontinuous functions
under the aw-distance. Even though epi-splines furnish a means to approximate
arbitrary semicontinuous functions using a finite number of parameters, they
suffer from the need to partition $\reals^n$ into a finite number of subsets.
In the present paper, we show that semicontinuous functions can be
approximated by piecewise affine functions that are defined {\it without
specifying a partition} and that are characterized structurally as being the
difference of two functions of the form $x\mapsto \max_{k=1, \dots, p}
\langle a^k, x\rangle + \alpha_k$. Consequently, we refer to these piecewise
affine functions as {\it mesh free}; the domain of each affine component
adapts and is not preselected. This is a significant feature for medium- and
high-dimensional problems, where representative low-dimensional subspaces
need to be discovered and exploited and standard polynomial approximations become challenging (see \cite{ZhangEtAl.15} for some progress in such directions). Our approximation result for usc
functions extends the well-known fact that every continuous function on a
convex compact set is the limit of a uniformly convergent sequence of dc
functions, which can be traced back to the local property of dc functions
established by \cite{Hartman.59}; see for example Proposition 2.3 in
\cite{HorstThoai.99}.

Covering numbers express the size of a class of functions in terms of the
smallest number of balls with a certain radius needed to cover the class and
are central to most consistency, rate of convergence, and error analysis in
(non)parametric estimation and machine learning; see for example
\cite{vanderVaartWellner.96,vandeGeer.00,GuoetAl.02}. The pioneering work
\cite{KolmogorovTikhomirov.61,BirmanSolomjak.67} deal with continuous and
smooth functions; see \cite{Pontil.03} for a more recent discussion.
Functions of bounded variation are considered in
\cite{BartlettKulkarniPosner.97} and analytic functions in \cite{Brudnyi.10}.
An upper estimate for the covering numbers of the unit ball of Gaussian
reproducing kernel Hilbert spaces is given in \cite{Zhou.02}, with further
refinements and applications in \cite{WangHuangLuoChen.09,Kuhn.11}. Covering
numbers of sets of convex functions are established in
\cite{Dudley.74,Bronshtein.76}, with significant improvements in
\cite{GuntuboyinaSen.13}. The present paper establishes an upper bound on the
covering numbers of bounded classes of usc functions under the aw-distance
and show that it is sharp within a logarithmic factor.

Although sample average approximations are often used to solving stochastic
optimization problems, it remains challenging to assess the quality of a
solution obtained through such approximations. Upper and lower bounds on
minimum values can be computed using the approaches in
\cite{HigleSen.96b,NorkingPflugRuszczynski.98,MMW.99,BayraksanMorton.06} (see
also \cite[Sect. 5.6]{ShapiroDentchevaRuszczynski.14}), at least when problem
relaxations can be solved to near global optimality. Validation approaches
based on optimality conditions are found in
\cite{HigleSen.91,SHo.98,Royset.12,LuLiuYinZhang.17,LammLu.18} and
\cite[Sect. 5.6]{ShapiroDentchevaRuszczynski.14}. Rates of convergence of
optimization problems with Lipschitz continuous objective functions defined
on a compact subset of $\reals^n$ are given in \cite[Sect.
5.3]{ShapiroDentchevaRuszczynski.14}. We leverage the results on covering
numbers to establish confidence regions of optimal solutions of
infinite-dimensional stochastic optimization problems defined on spaces of
usc functions without assuming Lipschitz continuity. The result is novel even
when specialized to finite dimensions. For a H\"{o}lder continuous case, we
obtain, in some sense, a stronger result.

After a section laying out notation and terminology, we proceed in Section 3
with the result on piecewise affine approximations and its applications.
Section 4 establishes bounds on covering numbers. Section 5 constructs
confidence regions and discusses rates of convergence for solutions of
stochastic optimization problems and their applications to nonparametric
estimation. The paper ends with an appendix supplementing a proof.

\section{Preliminaries}

In some applications, it would be natural and beneficial to consider usc
functions defined only on a strict subset of $\reals^n$ and their extensions
to the whole $\reals^n$ by assigning the value $-\infty$ may not be
meaningful. For example, if an usc function represents a necessarily
nonnegative probability density, then such an assignment would not imply a
useful extension. Consequently, we develop most results for usc functions
defined on a nonempty closed subset $S\subset\reals^n$, which could be all of
$\reals^n$, and is assumed to include the origin. Throughout, $S$ will be
such a set and the analysis will usually take place on the metric spaces
$(S,\|\cdot - \cdot\|_\infty)$ and $(S\times\reals,\|\cdot -
\cdot\|_\infty)$; the difference from the usual $(\reals^n,
\|\cdot-\cdot\|_\infty)$ is anyhow minor and will be highlighted when
significant. Of course, the sup-norm can be replaced by any other norm, but
this choice simplifies some expressions in Section 4. Likewise, the
assumption $0\in S$ can be relaxed, with the introduction of additional
notation better avoided here. The facts of this section can be found in or
deduced from \cite[Chapter 7]{VaAn} and \cite{Royset.18}.

We recall that $\hypo f = \{(x,\alpha)\in S\times \reals~|~f(x) \geq
\alpha\}\subset S\times \reals$ is the {\it hypograph} of a function
$f:S\to\Reals = [-\infty,\infty]$. The collection of {\it usc functions} on
$S$ is denoted by
\[
\uscfcns(S) = \{f:S\to \Reals~|~\hypo f \mbox{ is nonempty and closed}\}.
\]
We let $\nats = \{1, 2, \dots\}$. The {\it outer limit} of a sequence of sets
$\{A^\nu, \nu\in\nats\}$ in a topological space, denoted by $\nOutLim A^\nu$,
is the collection of points to which a subsequence of $\{a^\nu\in A^\nu,
\nu\in\nats\}$ converges. The {\it inner limit}, denoted by $\nInnLim A^\nu$,
is the collection of points to which a sequence $\{a^\nu\in A^\nu,
\nu\in\nats\}$ converges. If both limits exist and are equal to $A$, we say
that $\{A^\nu, \nu\in\nats\}$ {\it set-converges} to $A$ and write $A^\nu\to
A$ or $\nLim A^\nu = A$. We denote by $\nt A$ and $\cl A$ the interior and
closure of $A$, respectively.

For $f^\nu, f\in \uscfcns(S)$,
\[
f^\nu \mbox{ {\it hypo-converges} to } f, \mbox{ written } f^\nu\to f \Longleftrightarrow \hypo f^\nu \to \hypo f.
\]
Set-convergence of hypographs in this case, and therefore also
hypo-convergence, is equivalent to having
\begin{align}
  &\forall x^\nu\in S\to x, ~\nlimsup f^\nu(x^\nu) \leq f(x)\label{eqn:hypolimsup}\\
  &\forall x\in S, ~\exists x^\nu\in S\to x \mbox{ with } \nliminf f^\nu(x^\nu) \geq f(x).\label{eqn:hypoliminf}
\end{align}
The {\it Attouch-Wets (aw) distance} $\setd$, which quantifies the distance
between usc functions in terms of a distance between their hypographs,
metrizes hypo-convergence. Specifically, for $f,g\in \uscfcns(S)$, it is
defined as
\begin{equation*}
\setd(f,g) = \int_{0}^\infty \setd_\rho(f,g) e^{-\rho} d\rho,
\end{equation*}
where, for $\rho\geq 0$, the {\it $\rho$-aw-distance}
\begin{equation*}
\setd_\rho(f,g) = \nmax_{z\in \rho\ball_\infty} \big|\dist_\infty\big(z, \hypo f\big) - \dist_\infty\big(z, \hypo g\big)\big|,
\end{equation*}
with $\dist_\infty(z, A)$ being the usual point-to-set distance between a
point $z\in S\times \reals$ and a set $A\subset S\times \reals$ under the
sup-norm, $\rho\ball_\infty = \ball_\infty(0,\rho)$, with $\ball_\infty(\bar
z,\rho) = \{z \in S\times \reals~|~ \|\bar z - z\|_\infty \leq \rho\}$ for
any $\bar z\in S\times \reals$. Since the meaning will be clear from the
context, we also write $\ball_\infty(\bar x,\rho) = \{x \in S~|~ \|\bar x -
x\|_\infty \leq \rho\}$ with $\bar x\in S$. For any nonempty closed set
$S\subset\reals^n$, $(\uscfcns(S),\setd)$ is a complete separable metric
space. Every closed and bounded subset $F\subset (\uscfcns(S),\setd)$ is
compact. Moreover, for all $f,g\in \uscfcns(S)$,
\begin{equation}\label{eqn:distbds}
\big|\dist_\infty(0,\hypo f) - \dist_\infty(0,\hypo
g)\big| \leq \setd(f,g) \leq \max\big\{\dist_\infty(0,\hypo f), \dist_\infty(0,\hypo
g)\big\} + 1
\end{equation}
and, thus, if $f,g\geq 0$, then $0\leq\setd(f,g) \leq 1$. We also see that a
sufficient condition for $F$ to be bounded is that there exists
$(x,\alpha)\in S\times\reals$ such that $f(x)\geq \alpha$ for all $f\in F$.

If not specified otherwise, the index $\nu$ runs over $\nats$ so that
$x^\nu\to x$ means that the whole sequence $\{x^\nu, \nu\in\nats\}$ converges
to $x$. Let
\[
\cN_\infty^\grill \mbox{ be all the subsets of $\nats$ determined by subsequences},
\]
i.e., $N\in \cN_\infty^\grill$ is an infinite collection of strictly
increasing natural numbers. Thus, $\{x^\nu, \nu\in N\}$ is a subsequence of
$\{x^\nu, \nu\in\nats\}$; its convergence to $x$ is noted by $x^\nu\Nto x$.

A similar development is available for functions defined {\it on} the metric
space $(\uscfcns(S),\setd)$. However, we adopt a slightly different set-up
that highlights the role of domains of definition. For $F,F^\nu\subset
\uscfcns(S)$, the functions $\phi^\nu:F^\nu\to \Reals$ {\it epi-converge to }
$\phi:F\to\Reals$ whenever
\begin{align*}
  &\forall N\in \cN_\infty^\grill \mbox{ and }  f^\nu\in F^\nu\Nto f, ~\nliminf_{\nu\in N} \phi^\nu(f^\nu) \geq \phi(f) \mbox{ if } f\in F \mbox{ and } \phi^\nu(f^\nu)\Nto\infty \mbox{ otherwise}\\
  &\forall f\in F, ~\exists f^\nu\in F^\nu\to f \mbox{ with } \nlimsup \phi^\nu(f^\nu) \leq \phi(f).
\end{align*}
For $\epsilon\geq 0$, $\epsilon\mbox{-}\nargmin_{f\in F} \phi(f) = \{f\in
F~|~\phi(f) \leq \ninf_{g\in F} \phi(g) + \epsilon\}$, with the usual
extended real-valued calculus in play when needed. We deviate slightly from
the convention in \cite{VaAn} by setting $\epsilon\mbox{-}\nargmin_{f\in F}
\phi(f) = F$ when $\phi(f) = \infty$ for all $f\in F$. This is tenable
because we restrain from extending functions to the whole space and $\infty$
 is not assigned a special role in that regard. Consequently, we alleviate
 the need for checking that functions are finite at least somewhere, which in
 an infinite-dimensional setting may require excessively strong assumptions.

\section{Piecewise Affine Approximations}

In this section, we establish that every usc function on $S\subset\reals^n$
can be approximated by piecewise affine functions with a particular structure
under the additional assumption that $S$ is convex. For $\rho\in [0,\infty)$
and $q\in \nats$, let
\begin{align*}
&\pafcns^{q}_\rho(S) = \Big\{f:S\to [-\infty,\infty)~\Big|~\exists a^k,b^k\in\reals^n, \alpha_k,\beta_k\in \reals,~ k=1, \dots, q \mbox{ such that }\\
  & ~f(x) = \nmax_{k=1, \dots, q} \big[\langle a^k, x\rangle + \alpha_k\big] - \nmax_{k=1, \dots, q} \big[\langle b^k, x\rangle +
\beta_k\big]~\forall x\in S\cap \rho\ball_\infty;~ f(x) = -\infty \mbox{ otherwise}\}.
\end{align*}
A function in $\pafcns^{q}_\rho(S)$ is a difference of pointwise maxima of
affine functions on $S \cap \rho\ball_\infty$ and therefore is finite and
continuous on that set. We say that $U\subset\uscfcns(S)$ is {\it hypo-dense}
in $\uscfcns(S)$ if every $f\in \uscfcns(S)$ is the limit of a
hypo-converging sequence $\{f^\nu\in U, \nu\in\nats\}$.

\begin{theorem}{\rm (piecewise affine approximations).}\label{thm:pwa}
Suppose that $S$ is convex and $\rho^\nu\in [0,\infty)$ as well as $q^\nu\in
\nats$ tend to $\infty$. Then,
\[
\bigcup_{\nu\in \nats} \pafcns^{q^\nu}_{\rho^\nu}(S) \mbox{ is hypo-dense in } \uscfcns(S).
\]
\end{theorem}
\state Proof. Let $f\in \uscfcns(S)$. We construct a sequence in
$\cup_{\nu\in\nats} \pafcns^{q^\nu}_{\rho^\nu}(S)$ that hypo-converges to
$f$. Let $f^\nu:S\to \Reals$ have $f^\nu(x) = \min\{f(x),\nu\}$ for all $x\in
S$. Clearly, $f^\nu\to f$; recall that $\to$ always denotes hypo-convergence when written between usc functions. For any $\nu\in\nats$, $f^\nu$ is (upper)
prox-bounded so that for every $\lambda>0$, the (upper) Moreau-envelope
$e_\lambda f^\nu:S\to\reals$ of $f^\nu$, which is given by
\[
e_\lambda f^\nu(x) = \nsup_{y\in S} f^\nu(y) - \frac{1}{2\lambda} \|y-x\|_2^2,
\]
is finite and continuous. Moreover, $e_\lambda f^\nu\to f^\nu$ as
$\lambda\downto 0$ (see for example the discussion after Proposition 7.4 in \cite{VaAn}). Thus, there exists $\{\lambda^\nu>0, ~\nu\in\nats\}\to 0$
such that
\[
e_{\lambda^\nu} f^\nu \to f.
\]
Next, we define $\phi^\nu:S\to \Reals$ as
\[
\phi^\nu(x) = e_{\lambda^\nu} f^\nu(x) ~\forall x\in S\cap \rho^\nu\ball_\infty \mbox{ and } \phi^\nu(x) = -\infty \mbox{ otherwise}.
\]
Since $\rho^\nu\ball_\infty\to \reals^n$, we also have that $\phi^\nu\to f$.

Every real-valued continuous function on a convex compact subset of
$\reals^n$ is the limit in the sup-norm of finite-valued dc functions defined
on the same set; see for example \cite[Prop. 2.3]{HorstThoai.99}.
Consequently, for every $\nu\in\nats$, there exist convex functions
$\{g^\nu_\mu,h^\nu_\mu:S\to\Reals, ~\mu\in\nats\}$, finite on $S\cap
\rho^\nu\ball_\infty$, with the property that
\[
\nsup_{x\in S\cap \rho^\nu\ball_\infty} \big|\phi^\nu(x) - [g^\nu_\mu(x) - h^\nu_\mu(x)]\big| \to 0 \mbox{ as } \mu\to \infty.
\]
Let $\psi^\nu_\mu:S\to\Reals$ be defined by
\[
\psi^\nu_\mu(x) = g^\nu_\mu(x) - h^\nu_\mu(x) \mbox{ for } x\in S\cap \rho^\nu\ball_\infty \mbox{ and }  \psi^\nu_\mu(x) = - \infty \mbox{ otherwise.}
\]
Since already $\phi^\nu\to f$, we can construct $\{\mu^\nu \in \nats, ~\nu\in
\nats\}\to \infty$ as $\nu\to \infty$ such that $\psi^\nu_{\mu^\nu} \to f$.
Let $\psi^\nu = \psi^\nu_{\mu^\nu}$.

The convex functions $\{g^\nu_{\mu^\nu},h^\nu_{\mu^\nu}, ~\nu\in\nats\}$ are
lsc and proper. Let $g^\nu = g^\nu_{\mu^\nu}$ and $h^\nu = h^\nu_{\mu^\nu}$.
Consequently, for every $\nu\in\nats$,
\[
g^\nu(x) = \nsup_{(a,\alpha)\in A(g^\nu)} \{\langle a,x\rangle + \alpha\} \mbox{ and }   h^\nu(x) = \nsup_{(a,\alpha)\in A(h^\nu)} \{\langle a,x\rangle + \alpha\} \mbox{ for } x\in S,
\]
where for $u:S\to\Reals$,
\[
A(u) = \{(a,\alpha)\in \reals^n\times
\reals~|~\langle a,x\rangle + \alpha \leq u(x)~\forall x\in S\}.
\]
Since $A(g^\nu),A(h^\nu)\subset\reals^{n+1}$, which is separable, there exist increasing sets
$\{A^\mu(g^\nu), A^\mu(h^\nu), ~\mu\in\nats\}$, each with finite cardinality,
such that
\[
\bigcup_{\mu\in\nats} A^\mu(g^\nu) \mbox{ is dense in } A(g^\nu) \mbox{ and } \bigcup_{\mu\in\nats} A^\mu(h^\nu) \mbox{ is dense in } A(h^\nu).
\]
For $\nu, \mu\in\nats$, we define $\tilde g^\nu_\mu,\tilde h^\nu_\mu:S\to
\Reals$ by setting
\[
\tilde g_\mu^\nu(x) = \nmax_{(a,\alpha)\in A^\mu(g^\nu)} \{\langle a,x\rangle + \alpha\} \mbox{ and }   \tilde h_\mu^\nu(x) = \nmax_{(a,\alpha)\in A^\mu(h^\nu)} \{\langle a,x\rangle + \alpha\} \mbox{ for } x\in S\cap \rho^\nu\ball_\infty,
\]
and $\tilde g_\mu^\nu(x)=\tilde h_\mu^\nu(x)=\infty$ otherwise. The
characterization of hypo-convergence in \eqref{eqn:hypolimsup} and \eqref{eqn:hypoliminf} enables us to conclude that
for all $\nu\in\nats$,
\[
-\tilde g_\mu^\nu \to -g^\nu \mbox{ and } -\tilde h_\mu^\nu\to -h^\nu \mbox{ as } \mu\to \infty,
\]
with pointwise convergence holding as well on $S\cap \rho^\nu\ball_\infty$.

Let $\tilde\psi^\nu_\mu:S\to\Reals$ be defined by
\[
\tilde \psi^\nu_\mu(x) = \tilde g^\nu_\mu(x) - \tilde h^\nu_\mu(x) \mbox{ for } x\in S\cap \rho^\nu\ball_\infty \mbox{ and }  \tilde \psi^\nu_\mu(x) = - \infty \mbox{ otherwise.}
\]
If $\tilde \psi_\mu^\nu \to \psi^\nu$ as $\mu\to \infty$, then we can
construct $\{\mu^\nu\in\nats, ~\nu\in\nats\}\to \infty$ such that
\[
\tilde \psi^\nu_{\mu^\nu} \to f
\]
because already $\psi^\nu\to f$. Since $\tilde \psi^\nu_{\mu^\nu}\in
\pafcns^{q^\nu}_{\rho^\nu}(S)$, with $q^\nu$ being the largest cardinality of
$A^{\mu^\nu}(g^\nu)$ and of $A^{\mu^\nu}(h^\nu)$, the conclusion will follow.

It only remains to establish that $\tilde \psi_\mu^\nu\to \psi^\nu$ as
$\mu\to \infty$. For this purpose, we again leverage the characterization of hypo-convergence in \eqref{eqn:hypolimsup} and \eqref{eqn:hypoliminf}. Suppose that $x^\mu\in S\cap
\rho^\nu\ball_\infty\to x$. Then,
\begin{align*}
\nlimsup_\mu \big( \tilde g_\mu^\nu(x^\mu) - \tilde h_\mu^\nu(x^\mu) \big)& = \nlimsup_\mu \tilde g_\mu^\nu(x^\mu) - \nliminf_\mu \tilde h_\mu^\nu(x^\mu)\\
& \leq \nlimsup_\mu g^\nu(x^\mu) - h^\nu(x) \leq g^\nu(x) - h^\nu(x),
\end{align*}
where we use the facts that $\tilde g_\mu^\nu$ lower bounds $g^\nu$, $-\tilde
h_\mu^\nu\to -h^\nu$, and  $g^\nu$ is continuous on $S\cap
\rho^\nu\ball_\infty$. Also, for $x\in S\cap \rho^\nu\ball_\infty$,
\[
\nliminf_\mu \big( \tilde g_\mu^\nu(x) - \tilde h_\mu^\nu(x) \big) = \nliminf_\mu \tilde g_\mu^\nu(x) - \nlimsup_\mu \tilde h_\mu^\nu(x) \geq g^\nu(x) - h^\nu(x).
\]
These inequalities are trivially satisfied for sequences outside of $S\cap
\rho^\nu\ball_\infty$. Hence, the assertion is established.\eop

We illustrate the usefulness of the theorem in the solution of {\it  function
identification problems} of the form
\begin{equation*}\label{eqn:fid}
\mbox{(FIP)}~~~~\min_{f\in F} \phi(f), \mbox{ where } F\subset \uscfcns(S) \mbox{ and } \phi:F\to\Reals.
\end{equation*}
These problems arise in nonparametric estimation, spatial statistics, and
curve fitting; see \cite{RoysetWets.14} for applications to estimation of
financial curves, electricity demand, commodity prices, and uncertainty in
physical systems. For example, if $F$ is a class of $n$-dimensional
probability density functions and $\phi(f) = -\frac{1}{m}\sum_{j=1}^m \log
f(x^j)$, then any minimizer of (FIP) is a maximum likelihood estimate based
on the data $x^1, \dots x^m\in S$. When $\phi(f) = \frac{1}{m}\sum_{j=1}^m
(y^j - f(x^j))^2$, a minimizer furnishes a least-squares fit of the data
$\{(x^j, y^j)\in S\times \reals, ~j=1, \dots, m\}$ over the class $F$. We
refer to \cite{RoysetWets.15b,RoysetWets.18} for numerous examples.
Further unexplored applications for usc functions and their piecewise affine approximations may arise in stochastic and robust optimization where problems can be formulated over spaces of decision rules and be approximated using polynomial and piecewise polynomial functions \cite{BampouKuhn.12}, finite collections of policies \cite{HanasusantoKuhnWiesemann.15}, and linear decision rules, possibly in a higher dimensional spaces \cite{GeorghiouWiesemannKuhn.15}. In
special cases of (FIP), such as when $F$ consists of concave functions only, one might
be able to reformulate the problem as an equivalent finite-dimensional one; see \cite{CuleSamworthStewart.10a} for the context of maximum
likelihood estimation over the log-concave class and \cite{SeijoSen.11} for
least-squares regression over convex functions. However, this is not possible
in general and we need to settle for approximations.

For $\rho^\nu\in [0,\infty)$ and $q^\nu\in\nats$, we consider the {\it
approximating function identification problem}
\begin{equation*}\label{eqn:fidapprox}
\mbox{(FIP)}^\nu~~~~\min_{f\in F^\nu} \phi(f), \mbox{ where } F^\nu = F\cap \pafcns^{q^\nu}_{\rho^\nu}(S).
\end{equation*}
Every function in $\pafcns^{q^\nu}_{\rho^\nu}(S)$ is described by
$2q^\nu(n+1)$ parameters\footnote{We stress that $\nu$ is an index and not
the power of $q$.}. Thus, $\mbox{(FIP)}^\nu$ is equivalent to a
finite-dimensional optimization problem with the same number of variables.
The number grows only linearly in $n$, which makes the approach promising for
high-dimensional problems. In comparison, approximations based on epi-splines
(see \cite{RoysetWets.15b,RoysetWets.18}) require a preselected partition of
$S$ not easily decided on in a computationally tractable manner beyond four
or five dimensions. The piecewise approximations in $\mbox{(FIP)}^\nu$ are
mesh free, with the domain of each affine component adapting to the problem
at hand. Thus, we expect to be able to identify and leverage low-dimensional
structures, if present, when solving (FIP) by means of $\mbox{(FIP)}^\nu$.

It is apparent that an application may demand approximations also of the
objective function in (FIP), which we address in Section 5 for the central
case of stochastic optimization where $\phi = \Ex[\psi(\bfxi,\cdot)]$; the
expectation with respect to the distribution of a random element $\bfxi$ is
denoted by $\Ex$. Here, we concentrate on the application of Theorem
\ref{thm:pwa} to justify $\mbox{(FIP)}^\nu$.

Let $\phi^\nu:F^\nu\to \Reals$ be the function defined by $\phi^\nu(f) =
\phi(f)$ for $f\in F^\nu$. Suppose that $S$ is convex, $F$ is a nonempty and
{\it solid} subset of $(\uscfcns(S),\setd)$, i.e., $F=\cl (\nt F)$, and
$\phi:F\to\Reals$ is continuous on $F$. Then, a standard argument (see for
example the proof of Theorem 3.16 in \cite{RoysetWets.15b}) in conjunction
with Theorem \ref{thm:pwa} establishes that
\[
\phi^\nu \mbox{ epi-converges to } \phi \mbox{ provided that } \rho^\nu,q^\nu\to\infty.
\]
Thus, when $\{\epsilon^\nu\geq 0,~ \nu\in\nats\}\to 0$,
\[
\nOutLim \Big( \epsilon^\nu\mbox{-}\nargmin_{f\in F^\nu} \phi(f) \Big) \subset \nargmin_{f\in F} \phi(f),
\]
which can be deduced, for example, from \cite{Royset.18}. The constraint
qualification that $F$ is solid cannot be relaxed without introducing some
other assumption. Obviously, $F\cap \pafcns^{q^\nu}_{\rho^\nu}(S)$ can, in
general, be empty for all $\nu$, but when $F$ is solid this is ruled out.

In view of this discussion, the challenge of solving an infinite-dimensional
function identification problem from the broad class (FIP) is shifted to that
of obtaining a near-minimizer of a finite-dimensional problem. Of course, the
difficulty of that task depends on the specific properties of $\phi$ and $F$.
Typically, $\mbox{(FIP)}^\nu$ would be nonconvex, but the special affine
structure of functions in $\pafcns^{q^\nu}_{\rho^\nu}(S)$ is bound to be
important in developing computational procedures. Initial efforts in that
direction are already found in \cite{CuiPangSen.18}, which presents several
algorithms with guarantees to obtain at least certain stationary points and
numerical results from nonparametric least-squares regression, as well as in \cite{Miller.19}, which approximates functions in up to $n=41$ dimensions using piecewise affine functions in the context of nonparametric support vector machines. The nonconvexity of $\mbox{(FIP)}^\nu$ encountered in \cite{Miller.19} appears to be only moderately challenging and handled by common randomization strategies.  

\section{Covering Numbers}

It is well known that every bounded $F\subset(\uscfcns(S), \setd)$ has a
finite cover by virtue of being totally bounded. However, this is not
sufficient to establish certain rates of convergence results for sample
average approximations of stochastic optimization problems of the form
$\min_{f\in F} \Ex[\psi(\bfxi,f)]$. It is usually necessary to bound for all
$\epsilon>0$ the {\it covering number} of $F$, denoted by $N(F,\epsilon)$,
which is the smallest number of closed $\setd$-balls of radius $\epsilon$ needed to
cover $F$. We next provide such a bound and show that it is nearly sharp.
Section 5 applies the result to establish rates of convergence of minimizers
of stochastic optimization problems.

We start by recording a useful estimate of the hypo-distance. For $f,g\in
\uscfcns(S)$ and $\rho\geq 0$, we define the auxiliary quantity
\begin{align*}
  \hatsetd_\rho(f,g) = \ninf\Big\{ \tau\geq 0 ~\Big|~ &\nsup_{y\in\ball_\infty(x,\tau)} g(y) \geq \min\{f(x), \rho\} + \tau, \forall x\in \rho\ball_\infty \mbox{ with } f(x) \geq -\rho\\
                                                  &\nsup_{y\in\ball_\infty(x,\tau)} f(y) \geq \min\{g(x), \rho\} + \tau, \forall x\in \rho\ball_\infty \mbox{ with } g(x) \geq -\rho \Big\}.
\end{align*}
As the notation indicates, $\hatsetd_\rho$ is closely related to $\setd_\rho$
(see Proposition 3.1 in \cite{Royset.18}) and therefore also to $\setd$. We
record the relevant properties next.

\begin{lemma}\label{lemma:distestim}
  For $f,g\in \uscfcns(S)$ and $\rho\geq 0$,
  \[
e^{-\rho} \hatsetd_\rho(f,g) \leq \setd(f,g) \leq (1-e^{-\rho})\hatsetd_{2\rho+\delta}(f,g) + e^{-\rho}(\delta+\rho+1),
  \]
where $\delta = \max\{\dist_\infty(0,\hypo f), \dist_\infty(0,\hypo g)\}$.
\end{lemma}
\state Proof. The results can be deduced from Propositions 3.1 and 3.2 in
\cite{Royset.18}.\eop

As mentioned in Section 1, epi-splines \cite{RoysetWets.15b,Royset.18}
furnish a dense subset of classes of semicontinuous functions and associated
error bounds are known. We leverage these results here. Although the piecewise affine functions of Section 3 are also dense in the usc functions, they have unknown error and cannot presently serve as the basis for the construction in the proof of the next theorem. This is anyhow less critical as we see through a lower bound result (Theorem \ref{thm:mainlow} below) that the obtained upper bound on covering numbers is within a logarithmic factor of being sharp.

For any
$f:S\to\Reals$ and $x\in S$, let $\nliminf_{\bar x\to x} f(\bar x) =
\nlim_{\delta\downarrow 0}\ninf_{\bar x\in\ball_\infty(x,\delta)} f(\bar x)$.
Epi-splines are defined in terms a finite collection of subsets of $S$. A
finite collection $R_1, R_2, ..., R_K$ of open subsets\footnote{Recall that
``open'' here is according to the metric space $(S,\|\cdot-\cdot\|_\infty).$}
of $S$ is a {\it partition} of $S$ if $\cup_{k=1}^K \cl R_k = S$ and $R_k
\cap R_l = \emptyset$ for all $k\neq l$. Specifically, an {\it epi-spline}
$s:S\to\reals$, with partition $\cR=\{R_1, \dots, R_K\}$ of $S$, is a
function that
\begin{align*}
&\mbox{on each $R_k$, $k=1, ..., K$, takes a constant real number as value,}\\
&\mbox{and for every $x\in S$,} \mbox{ has } s(x) = \nliminf_{x'\to x} s(x').
\end{align*}
The family of all such epi-splines is denoted by $\espl(\cR)$. Epi-splines
are lsc by construction and approximate lsc functions in the sense of
epi-convergence. Since the present setting involves usc functions and
hypo-convergence, we ``reorientation'' and introduce minus in some
expressions. We refer to \cite{RoysetWets.15b,Royset.18} for further
information and extensions that go beyond these zeroth order epi-splines and
also beyond $\reals^n$.

\begin{proposition}\label{prop:rateofepiconv} For a partition $\cR =\{R_1, \dots, R_K\}$ of  $S$ and
$\rho\geq 0$, we have that for every $f\in \uscfcns(S)$, there exists an
$s\in \espl(\cR)$ such that
  \begin{equation*}
     \hatsetd_\rho(f,-s) \leq \mu_\rho(\cR)  = \ninf\big\{ \tau\geq 0~|~ R_k \subset \ball_\infty(x,\tau) \mbox{ for all } x\in  \rho\ball_\infty \mbox{ and } k \mbox{ satisfying } x \in \cl R_k\big\}.
  \end{equation*}
If $\mu_\rho(\cR) \leq \rho$, then $s$ can be taken to satisfy $-\rho' \leq
s(x) \leq \max\{-\rho', \min[\rho', -f(x)]\}$ for any $\rho'>\rho$ and $x\in
S$.
\end{proposition}
\state Proof. The first part of the proposition is a direct application of
\cite[Theorem 5.9]{Royset.18}. The fact that $s$ can be taken to satisfy
$-\rho' \leq s(x) \leq \max\{-\rho',\min[\rho',-f(x)]\}$ for any $\rho'>\rho$
follows from an examination of that theorem's proof.\eop

The quantity $\mu_\rho(\cR)$ is the {\it meshsize} of $\cR =\{R_1, \dots,
R_K\}$ and, essentially, quantifies the size of the largest $R_k$.

\begin{theorem}\label{thm:entropy}{\rm (covering numbers).}
For every bounded subset $F$ of $(\uscfcns(S), \setd)$, there exist $c\geq 0$
and $\bar \epsilon>0$ (both independent of $n$, the dimension of $S$) such
that
\begin{equation*}
  \log N(F,\epsilon) \leq \left(\frac{c}{\epsilon}\right)^n\left(\log \frac{1}{\epsilon}\right)^{n+1}  \mbox{ for all } \epsilon \in (0, \bar \epsilon].
\end{equation*}
\end{theorem}
\state Proof. Since $F$ is bounded, there exists an $r>0$ such that
$\dist_\infty(0,\hypo f) \leq r$ for all $f\in F$. Let
$\gamma_1,\gamma_2,\gamma_3
> 0$ be such that $\gamma_1+\gamma_2+\gamma_3 = 1$. Set $\bar \epsilon \in
(0,1)$ such that
\begin{equation}\label{eqn:covering1}
  \frac{2(r + 1)}{r}\left[\log \frac{1}{\epsilon} + \log \frac{1}{\gamma_1} +  \frac{r}{2} + \log\left(r + 1\right) \right] - 1>\gamma_2 \epsilon \mbox{ for all } \epsilon \in (0,\bar \epsilon].
\end{equation}
Fix $\epsilon \in (0,\bar \epsilon]$ and define $\rho$ to be the expression
on the left-hand side of \eqref{eqn:covering1}. We next construct a
partition of $S$ and set $\omega >1$ and
\begin{equation*}
  \nu = \left\lceil  \frac{2\omega\rho}{\gamma_2 \epsilon} \right\rceil,
\end{equation*}
where $\lceil a \rceil$ is the smallest integer no smaller than $a$. The
partition is obtained by dividing the box $[-\omega\rho,\omega\rho]^n\subset
\reals^n$ into $\nu^n$ boxes of equal size and then intersecting with $S$.
Let $K = \nu^n + 1$. Specifically, for $k=1, 2, \dots , \nu^n$, set
\[
R_k = \nt \Big(S\cap \prod_{i=1}^n (l_i^k, u_i^k)\Big), \mbox{ with } l_i^k =
2(k-1)\omega\rho/\nu - \omega\rho \mbox{ and } u_i^k  = l_i^k + 2\omega\rho/\nu
\]
so that $\cup_{k=1}^{K-1} \cl R_k = S\cap [-\omega\rho,\omega\rho]^n$. Also,
$R_K = \nt ( S\setminus [-\omega\rho,\omega\rho]^n)$. Again, we recall that
the interior and closure are taken relative to $(S,\|\cdot-\cdot\|_\infty)$.
We denote by $\cR = \{R_1, \dots, R_K\}$ this partition. Clearly,
$\mu_{\rho}(\cR) = 2\omega\rho/\nu$. Next, we consider a discretization of
parts of the range of functions and set
\begin{equation*}
  m = \left\lceil \frac{\omega\rho}{\gamma_3 \epsilon } \right\rceil + 1.
\end{equation*}
The points $\sigma_j = -\omega\rho + 2(j-1)\omega\rho/(m-1)$, $j=1, 2, \dots
, m$, discretize the interval $[-\omega\rho, \omega\rho]$. Let $F_0$ be the
collection of piecewise constant functions on $\cR$ with values in
$\{\sigma_1, \dots . \sigma_m\}$ defined as follows. If $f\in F_0$, then for
every $k\in \{1, \dots , K\}$ there exists a $j_k\in \{1, \dots , m\}$ such
that $f(x) = \sigma_{j_k}$ for $x\in R_k$ and $f(x) = \nlim_{\delta\downarrow
0} \nsup_{y\in \ball_\infty(x,\delta)} f(y)$ otherwise. By construction, $f$
is usc. Obviously, $F_0$ contains $m^K$ functions.
 We now show that every $f\in F$ has $\setd(f,f_0) \leq \epsilon$ for some
 $f_0\in F_0$.

Let $f\in F$ be arbitrary. By Proposition \ref{prop:rateofepiconv} and the
fact that $\mu_\rho(\cR) = 2\omega\rho/\nu \leq \gamma_2\epsilon < \rho$,
there exists $s\in \espl(\cR)$ such that
\begin{equation*}
  \hatsetd_{\rho}(f,-s) \leq \mu_{\rho}(\cR) \mbox{ and } -\omega\rho \leq s(x) \leq \max\{-\omega\rho, \min[\omega\rho,-f(x)]\} \mbox{ for } x\in S.
\end{equation*}
Since $\dist(0,\hypo f) \leq r$, there exists $x\in r\ball_\infty$ such that
$f(x)\geq -r$. Consequently, $-s(x) \geq \min\{\omega\rho,
\max[-\omega\rho,f(x)]\} \geq -r$. So we also have that $\dist_\infty(0,\hypo
-s) \leq r$.

Since $\epsilon,\gamma_1 \leq 1$,
\begin{equation*}
  \rho \geq \frac{2(r+1)}{r}\left[\frac{r}{2} + \log (r+1) \right] - 1 = r + \frac{2(r+1)}{r}\log(r+1) \geq r.
\end{equation*}
Thus, using the notation $\bar\rho = (\rho-r)/2$, Lemma \ref{lemma:distestim}
gives that
\[
  \setd(f,-s) \leq \hatsetd_{\rho}(f,-s) + e^{-\bar\rho}(r + \bar\rho + 1)
             \leq \mu_{\rho}(\cR) + e^{-\bar\rho}(r + \bar\rho + 1)
              =    2\omega\rho/\nu + e^{-\bar\rho}(r + \bar\rho + 1).
\]
In view of \cite[Theorem 3.17]{RoysetWets.15b}, there exists $f_0 \in F_0$
such that $\setd(-s,f_0) \leq \omega\rho/(m-1)$ since we can select $f_0$
such that $|s(x)-f_0(x)| \leq \omega\rho/(m-1)$ for all $x\in S$. The
triangle inequality then gives that
\begin{equation}\label{eqn:covering2}
  \setd(f,f_0) \leq \omega\rho/(m-1) + 2\omega\rho/\nu + e^{-\bar\rho}(r + \bar\rho + 1).
\end{equation}
It remains to show that the right-hand side is no greater than $\epsilon$. We
start with the last term in \eqref{eqn:covering2}. By concavity of the log-function, we have that
\begin{equation*}
  \log\left(\frac{1}{2}(\rho+r)+1\right) \leq \log\left(r+1\right) + \frac{\rho-r}{2r+ 2}.
\end{equation*}
Consequently,
\begin{align*}
  \log\left[e^{-\bar\rho}(r + \bar\rho + 1)\right] &=  \frac{1}{2}(r - \rho) + \log\left(\frac{1}{2}(\rho+r)+1\right) \leq \frac{1}{2}(r - \rho) + \log\left(r+1\right) + \frac{\rho-r}{2r+ 2}\\
 & = \frac{r}{2}-\frac{r(\rho+1)}{2(r+1)} + \log(r+1) = \log \gamma_1 \epsilon,
\end{align*}
where the last equality follows from inserting the expression for $\rho$.
Thus, $e^{-\bar\rho}(r + \bar\rho + 1)\leq \gamma_1\epsilon$. We next examine
the second term on the right-hand side of \eqref{eqn:covering2}. Inserting the expression for $\nu$, we obtain that
\begin{equation*}
  \frac{2\omega\rho}{\nu} \leq \gamma_2 \epsilon.
\end{equation*}
Finally, we consider the first term on the right-hand side of \eqref{eqn:covering2}. In view of the definition of $m$, we
have that
\begin{equation*}
  \frac{\omega\rho}{m-1} \leq \gamma_3 \epsilon.
\end{equation*}
Thus, $\setd(f,f_0) \leq \epsilon$ and we have established that $\setd$-balls
with radius $\epsilon$ and centered at points in $F_0$ cover $F$. The
logarithm of the number of functions in $F_0$ is $(\nu^n + 1)\log m$. At this
point, the order of the result is immediate. A possible expression for the
constant $c$ is obtained as follows. Let $c_1 = 2(r+1)/r$ and
\begin{equation*}
c_2 =   \frac{2(r + 1)}{r}\left[\log \frac{1}{\gamma_1} +  \frac{r}{2} + \log\left(r + 1\right) \right] - 1.
\end{equation*}
Thus, $\rho = c_1\log \epsilon^{-1} + c_2$. Moreover, let $c_3 =
2\omega/\gamma_2$ and $c_4 = \omega/\gamma_3$. Using these expressions, we
find that
\begin{equation*}
  (\nu^n + 1)\log m \leq \left[\left( c_1c_3 + \frac{c_2c_3+1}{\log \bar \epsilon^{-1}}  \right)^n \left(\frac{1}{\epsilon} \log \frac{1}{\epsilon} \right)^n   +1 \right] \log\left[ \left( c_1c_4 + \frac{c_2c_4+2}{\log \bar \epsilon^{-1}} \right) \frac{1}{\epsilon}\log \frac{1}{\epsilon} \right].
\end{equation*}
Let
\begin{equation*}
  c_5 = c_1c_3 + \frac{c_2c_3+1}{\log \bar \epsilon^{-1}} \mbox{ and } c_6 = c_1c_4 + \frac{c_2c_4+2}{\log \bar \epsilon^{-1}}.
\end{equation*}
We then find that
\begin{equation*}
  (\nu^n + 1)\log m \leq c_7^n \left(\frac{1}{\epsilon} \log \frac{1}{\epsilon} \right)^n \left[ \log c_6 + \log \frac{1}{\epsilon} + \log\log\frac{1}{\epsilon}  \right], \mbox{ where } c_7 = c_5 + \frac{1}{\bar \epsilon^{-1}\log \bar \epsilon^{-1}}.
\end{equation*}
Using the fact that $\log\log \epsilon^{-1}/\log \epsilon^{-1} \leq e^{-1}$
for $\epsilon \in (0,1)$, we obtain
\begin{equation}\label{eqn:covering3}
  (\nu^n + 1)\log m \leq c_7^n\left[\frac{\log c_6}{\log \bar \epsilon^{-1}} + 1 + e^{-1} \right] \frac{1}{\epsilon^n} \left(\log \frac{1}{\epsilon} \right)^{n+1},
\end{equation}
which gives a particular expression for $c$ in the theorem statement. Since
the choice of $\bar\epsilon$ is independent of $n$, this $c$ is independent
of $n$. For example, for $\bar\epsilon = 0.01$, $\omega = 0.00000001$, $r =
3.22$, then $c_7 = 13.5$ and the term in brackets in \eqref{eqn:covering3} evaluates to $2.3$.\eop

Although a comparison to the classical result of $O(\epsilon^{-n})$ for
Lipschitz continuous functions on bounded subsets, which goes back to
\cite{KolmogorovTikhomirov.61} (see for example \cite[Theorem
2.7.1]{vanderVaartWellner.96}), is not entirely relevant due to the different
settings, we note that our bound is only slightly worse (a logarithmic term)
for larger families of usc functions. We do not require any bound on the
variation of the functions and allow functions defined on all of $\reals^n$,
possibly extended real-valued. Still, the entropy integral\footnote{For the
significance of entropy integrals we refer to \cite{vanderVaartWellner.96}.}
$\int_0^{\bar \epsilon} \sqrt{\log N(F,\epsilon)} d\epsilon$ is finite only
for $n = 1$ and, therefore, these families are ``large,'' and increasingly so
as $n$ grows.

\begin{theorem}\label{thm:mainlow}{\rm (covering numbers; lower bound).}
For every $n\in\nats$, there exist a bounded subset $F\subset
\uscfcns(\reals^n)$ and corresponding $c\geq 0$ and $\bar \epsilon>0$
(independent of $n$) such that
\begin{equation*}
  \log N(F,\epsilon) \geq \left(\frac{c}{\epsilon}\right)^n \log \frac{1}{\epsilon} \mbox{ for all } \epsilon \in (0, \bar \epsilon].
\end{equation*}
\end{theorem}
\state Proof. See the appendix.\eop

In comparison with the upper bound of Theorem \ref{thm:entropy}, we see that
the lower bound differs by a logarithmic factor only and therefore the upper
bound is nearly sharp. The size of the bounded set $F$ in Theorem
\ref{thm:mainlow} does not have to be large. In fact, an examination of the
proof reveals that $F$ might be selected to have $\setd(0,f)\leq r$ for all
$f\in F$, with $r$ being only slightly above one. Here, $0$ is the function
in $\uscfcns(\reals^n)$ that is identical to zero.

The proof of Theorem \ref{thm:mainlow} constructs a collection of functions
which is finite on a grid of points in $[0,\rho]^n$, with $\rho>0$ and grid
points spaced roughly $\epsilon$ apart. At each of these grid points, a
function takes on a value among a set of discretized values between $-\rho$
and $0$, again spaced roughly $\epsilon$ apart. Outside these grid points,
the functions are assigned $-\infty$. It is clear that the number of such
functions is $(\rho/\epsilon)^\nu$, where $\nu = (\rho/\epsilon)^n$. Thus,
its logarithm is of the order $O(\epsilon^{-n} \log \epsilon^{-1})$.

\section{Applications to Stochastic Optimization and Statistical Estimation}

Suppose that  $(\Xi, \cA, \mathbb{P})$ is a complete probability space,
$F\subset \uscfcns(S)$ is closed, and $\psi:\Xi\times F\to \Reals$ is a
function with suitable properties as discussed below. We denote by boldface,
for example $\bfxi$, random elements with values in $\Xi$. Then, a {\it
function identification problem under uncertainty} takes the form
\[
\mbox{(FIP-U)}~~~~\min_{f\in F} \Ex[\psi(\bfxi,f)] = \int \psi(\xi,f) d\mathbb{P}(\xi).
\]
Section 3 furnishes some instances of $\psi$ in the context of probability
density estimation and regression, see also below, but there are numerous
other examples.

A {\it sample average approximation} of the problem leverages a sample
$\bfxi^1, \bfxi^2, \dots$ of independent random elements, each with values in
$\Xi$ and distributed according to $\mathbb{P}$, and leads to the {\it
approximating problem}
\[
\mbox{(FIP-U)}^\nu~~~~ \min_{f\in F} \frac{1}{\nu} \sum_{j=1}^\nu \psi(\bfxi^j, f).
\]
Under mild assumptions (see Proposition \ref{prop:consistency} below),
minimizers of $\mbox{(FIP-U)}^\nu$ tend to those of (FIP-U) almost surely.
However, a main challenge in the justification of such an approach is to
quantify the {\it rate} with which the error in solutions of
$\mbox{(FIP-U)}^\nu$ vanishes as $\nu$ grows. Before stating the results that rely on the covering numbers of the previous section, we
formalize the setting. The following definitions and facts are well known;
see for example\footnote{This reference states results only for finite
dimensions, but since $(F,\setd)$ is a complete separable metric space, with
compact balls, the proofs of the required results carry over nearly
verbatim.} \cite[Ch. 14]{VaAn}.

We say that $\psi:\Xi\times F\to \Reals$ is a {\it random lsc function} if
for all $\xi\in \Xi$, $\psi(\xi,\cdot)$ is lsc as a function on the metric
space $(F,\setd)$ and $\psi$ is measurable with respect to the product
sigma-algebra\footnote{On $(F,\setd)$, we adopt the Borel sigma-algebra.} on
$\Xi\times F$. A random lsc function $\psi:\Xi\times F\to \Reals$ is {\it
locally inf-integrable} if\footnote{For measurable $h:\Xi\to \Reals$, $\int
h(\xi)d\mathbb{P}(\xi) = \int \max\{0,h(\xi)\}d\mathbb{P}(\xi) - \int
\max\{0, -h(\xi)\}d\mathbb{P}(\xi)$, with $\infty-\infty = \infty$.}
\begin{equation*}
  \forall f\in F~\exists \rho>0 \mbox{ such that } \int \ninf_{g\in F} \{\psi(\xi,g)~|~\setd(g,f)\leq \rho\}~ d\mathbb{P}(\xi) > -\infty.
\end{equation*}
If $\psi:\Xi\times F\to \Reals$ is a locally inf-integrable random lsc
function, then $f\mapsto \Ex[\psi(\bfxi,f)]$ is well-defined, always greater
than $-\infty$, and lsc.

\begin{proposition}\label{prop:consistency} Suppose that  $(\Xi, \cA, \mathbb{P})$ is a complete probability space,
$F\subset \uscfcns(S)$ is closed, and $\psi:\Xi\times F\to \Reals$ is an
inf-integrable random lsc function. If $\bfxi^1, \bfxi^2, \dots$ is a
sequence of independent  random elements each with values in $\Xi$ and
distribution $\mathbb{P}$ and $\{\epsilon^\nu\geq 0, \nu\in\nats\}\to 0$,
then, almost surely,
\[
\nOutLim \Big(\epsilon^\nu\mbox{-}\nargmin_{f\in F} \frac{1}{\nu}\sum_{j=1}^\nu \psi(\bfxi^j,f)\Big)     \subset \nargmin_{f\in F} \Ex[\psi(\bfxi^1,f)].
\]
Moreover, if $F$ is bounded, then almost surely
\[
\lim \Big(\ninf_{f\in F} \frac{1}{\nu}\sum_{j=1}^\nu \psi(\bfxi^j,f)\Big) = \ninf_{f\in F} \Ex[\psi(\bfxi^1,f)] > -\infty.
\]
\end{proposition}
\state Proof. This is a consequence of a law of large numbers for lsc
functions and epi-convergence; see for example Proposition 7.1 in
\cite{RoysetWets.18}. \eop

In the following, we assume that $F\subset\uscfcns(S)$ is closed and bounded
as it results in some simplifications. In particular, $(F,\setd)$ then
becomes a {\it compact} metric space. The assumption is anyhow minor as it is
often acceptable in applications to impose a lower bound on the functions in
$\uscfcns(S)$ under considerations; see the remark in conjunction with
(\ref{eqn:distbds}).

The {\it excess} of a set $F_1\subset F$ over a set $F_2\subset F$ is given
by
\begin{equation*}
  \exs(F_1,F_2) = \nsup_{f\in F_1} \dist(f, F_2) \mbox{ if } F_1,F_2 \mbox{ are nonempty},
\end{equation*}
$\exs(F_1,F_2) = \infty$ if $F_1$ nonempty and $F_2$ empty, and
$\exs(F_1,F_2) = 0$ otherwise. Here, $\dist(f, F_2) = \ninf_{g\in F_2}
\setd(f,g)$ is the usual point-to-set distance in $(F,\setd)$. The
Pompeiu-Hausdorff distance $\mathbb{H}(F_1,F_2) = \max\{\exs(F_1,F_2),
\exs(F_2,F_1)\}$. Let the {\it level sets} of any $\phi:F\to \Reals$ be
denoted by
\[
\nlev_{\leq \delta} \phi = \{f\in F~|~\phi(f) \leq \delta\}.
\]
If $\psi:\Xi\times F\to \Reals$ is a random lsc function and $F_1,F_2\subset
F$ are closed, then
\[
\xi\mapsto\exs\big(\epsilon\mbox{-}\nargmin_{f\in F_1} \psi(\xi,f), ~F_2\big) \mbox{ and } \xi\mapsto\exs\big(F_2, ~\nlev_{\leq \delta} \psi(\xi,\cdot)\big)
\]
are random variables on $(\Xi, \cA, \mathbb{P})$ for any $\epsilon\geq 0$ and
$\delta\in\reals$.

When considering sample average approximations with sample size $\nu$, the
relevant probability space is the $\nu$-fold product space formed by
$(\Xi,\cA, \mathbb{P})$; the above definitions apply also for this
probability space. The {\it sample average function}
\[
\big((\xi^1, \dots , \xi^\nu),f\big)\mapsto \frac{1}{\nu}\sum_{j=1}^\nu \psi(\xi^j,f)
\]
is then a random lsc function on the product probability space provided that
$\psi>-\infty$ is a random lsc function on $(\Xi,\cA,\mathbb{P})$. Since this
is the case below, the following probabilistic statements are meaningful. The
measure on the product probability space is denoted by $P^\nu$ and the sample
space by $\Xi^\nu$.

We also need a quantitative result about differences between minimizers and
related quantities. The following result improves on \cite[Thms. 4.3 and
4.5]{Royset.18}; see also \cite{DevolderGlineurNesterov.10} for related results in the convex setting. We denote by $\ball(f,\rho) = \{g\in
\uscfcns(S)~|~\setd(f,g)\leq \rho\}$.

\begin{proposition}\label{prop:argmin} For a closed and bounded $F\subset\uscfcns(S)$, let $F_1,F_2\subset F$ be nonempty and
$\phi_1:F_1\to (-\infty,\infty]$ as well as $\phi_2:F_2\to (-\infty, \infty]$
be lsc functions on the metric space $(F,\setd)$. Suppose that for some
$\tau,\gamma\geq 0$,
\[
\ball(g,\gamma)\cap F_1\neq\emptyset \mbox{ and }  \ninf_{f\in\ball(g,\gamma)\cap F_1} \phi_1(f) \leq \phi_2(g) + \tau ~\forall g\in F_2.
\]
Then, for any $\delta\in\reals$,
\begin{equation}\label{eqn:levelerror}
\exs\big(\nlev_{\leq \delta} \phi_2, \nlev_{\leq \delta + \tau} \phi_1\big)
\leq \gamma.
\end{equation}
If in addition
\[
\ball(g,\gamma)\cap F_2\neq\emptyset \mbox{ and }  \ninf_{f\in\ball(g,\gamma)\cap F_2} \phi_2(f) \leq \phi_1(g) + \tau ~\forall g\in F_1,
\]
then for any $\epsilon\geq 0$,
\begin{equation}\label{eqn:argminerror}
    \exs\big(\epsilon\textrm{-}\nargmin_{f\in F_1} \phi_1(f), (\epsilon+2\tau)\textrm{-}\nargmin_{f\in F_2} \phi_2(f)\big) \leq \gamma.
  \end{equation}
\end{proposition}
\state Proof. Let $g\in \nlev_{\leq \delta} \phi_2$. Since $\phi_1$ is lsc
and $\ball(g,\gamma)$ is compact, there exists $f^\star \in
\ball(g,\gamma)\cap F_1$ such that
\[
\phi_1(f^\star)=\ninf_{f\in\ball(g,\gamma)\cap F_1} \phi_1(f) \leq \phi_2(g)
+ \tau \leq \delta+\tau.
\]
We have established that $f^\star \in \nlev_{\leq\delta + \tau} \phi_1$.
Thus, $\dist(g,\nlev_{\leq\delta + \tau} \phi_1)\leq \gamma$ and \eqref{eqn:levelerror} follows.

For \eqref{eqn:argminerror}, we note that there exists $f^\star \in \nargmin_{f\in
F_2} \phi_2(f)$ because $F_2$ is totally bounded. Thus,
\[
\ninf_{f\in F_1} \phi_1(f) \leq \ninf_{f\in \ball(f^\star,\gamma)\cap F_1} \phi_1(f) \leq \phi_2(f^\star) + \tau = \ninf_{f\in F_2} \phi_2(f) + \tau.
\]
Suppose that $g\in\epsilon\textrm{-}\nargmin_{f\in F_1} \phi_1(f)$. Again,
there exists $f^{\star\star} \in \ball(g,\gamma)\cap F_2$ such that
\[
\phi_2(f^{\star\star})=\ninf_{f\in\ball(g,\gamma)\cap F_2} \phi_2(f) \leq \phi_1(g) + \tau \leq \ninf_{f\in F_1} \phi_1(g) + \epsilon + \tau
\leq \ninf_{f\in F_2} \phi_2(g) + \epsilon +2\tau.
\]
We have established that $f^{\star\star}
\in(\epsilon+2\tau)\mbox{-}\nargmin_{f\in F_2}  \phi_2(f)$. Thus,
$\dist(g,(\epsilon+2\tau)\mbox{-}\nargmin_{f\in F_2} \phi_2(f))\leq \gamma$
and \eqref{eqn:argminerror} follows.\eop

We observe that if $\phi_1$ and $\phi_2$ in the proposition are pointwise within $\delta$ of each other uniformly on $F$, then $\tau$ can be set to $\delta$ and $\gamma$ to zero. However, the focus on uniform bounds is limiting as it rules out discontinuous functions and especially cases with $F_1 \neq F_2$.

\subsection{Confidence Regions}

We are then in a position to state the first of the two main results in this
section.

\begin{theorem}{\rm (confidence region).}\label{thm:confidence} For a complete probability space $(\Xi,\cA,\mathbb{P})$ and a closed and bounded set $F\subset\uscfcns(S)$,
suppose that $\psi:\Xi\times F\to (-\infty,\infty]$ is an inf-integrable
random lsc function, $\bfxi^1, \bfxi^2, \dots$ are independent random
elements, each with values in $\Xi$ and distributed according to
$\mathbb{P}$, and $\psi(\bfxi^1,f)$ is sub-exponential\footnote{A random
variable $Y$ is sub-exponential if for some $\lambda\geq 0$, $\Ex[\exp(\tau
(Y-\Ex Y))] \leq \exp(\tau^2\lambda^2/2)$ for all $|\tau|\leq 1/\lambda$. Another assumption that ensures a Bernstein-type large-deviation result could have been substituted here.}
for all $f\in F$. Given $\alpha\in (0,1)$ and $\delta > \ninf_{f\in F}
\Ex[\psi(\bfxi^1,f)]$, there exist $\bar \nu\in\nats$ and $c\in [0,\infty)$
such that for all $\nu\geq \bar \nu$
\[
P^\nu\bigg[\exs\Big(\nargmin_{f\in F} \Ex[\psi(\bfxi^1,f)], ~\nlev_{\leq \delta} \Big\{\frac{1}{\nu}\sum_{j=1}^\nu \psi(\bfxi^j, \cdot)\Big\}\Big)
\leq \frac{c(\log \nu)^{1+1/n}}{\nu^{1/n}}\bigg] \geq 1 - \alpha.
\]
\end{theorem}
\state Proof. Let $\phi:F\to \reals$
  have values $\phi(f) =
\Ex[\psi(\bfxi^1,f)]$, which is well-defined, lsc, and indeed finite valued
due the sub-exponential assumption. Let $\phi^\nu:\Xi^\nu\times F\to
(-\infty, \infty]$ have values $\phi^\nu((\xi^1, \dots , \xi^\nu), f) =
\nu^{-1}\sum_{j=1}^\nu \psi(\xi^j,f)$, which then is a random lsc function on
the product probability space. At a given $f\in F$, with probability one,
$\phi^\nu((\bfxi^1, \dots , \bfxi^\nu), f) < \infty$ because otherwise
$\Ex[\psi(\bfxi^1,f)]$ would not have been finite.

As in Theorem \ref{thm:entropy}, there is a finite number $N = N(F,\gamma_1)$
of closed balls in $(\uscfcns(S),\setd)$ with radius $\gamma_1>0$ and center
$f_k\in\uscfcns(S)$ covering $F$. Without loss of generality, we can assume
that $\ball(f_k,\gamma_1)\cap F \neq \emptyset$. Moreover, let $f_k^\star
\in\nargmin_{f\in \ball(f_k,\gamma_1)\cap F} \phi(f)$. Since
$\psi(\bfxi^1,f_k^\star)$ is sub-exponential, there exists by Bernstein's
inequality $\gamma_2 \in (0, \delta - \ninf_{f\in F} \phi(f))$ and $c_0>0$
such that
\[
P^\nu\big(|\phi^\nu((\bfxi^1, \dots , \bfxi^\nu), f_k^\star) - \phi(f_k^\star)| \geq \gamma_2\big) \leq 2e^{-\nu c_0\gamma_2^2} \mbox{ for all } k=1, \dots , N.
\]
Consequently, as long as
\begin{equation}\label{eqn:samplesize}
2Ne^{-\nu c_0\gamma_2^2} \leq \alpha~
\mbox{ or, equivalently, } ~\nu\geq \frac{\log N - \log(\alpha/2)}{c_0 \gamma_2^2}
\end{equation}
we have that
\[
P^\nu\bigg(\max_{k=1, \dots ,N} \Big|\phi^\nu((\bfxi^1,\dots , \bfxi^\nu), f_k^\star) - \phi(f_k^\star)\Big| \geq \gamma_2\bigg) \leq \alpha.
\]
Suppose that we have an event $(\xi^1, \dots , \xi^\nu)\in \Xi^\nu$ where
\[
\max_{k=1, \dots ,N} \Big|\phi^\nu((\xi^1,\dots , \xi^\nu), f_k^\star) - \phi(f_k^\star)\Big| < \gamma_2.
\]

Next, we apply Proposition \ref{prop:argmin} and start by establishing the
required condition. Let $g\in F$ and $\delta_0 = \ninf_{f\in F} \phi(f)$,
which is finite because $F$ is compact. Then, there exists $k^\star\in \{1,
\dots , N\}$ such that $g\in \ball(f_{k^\star},\gamma_1)$ and
\[
  \ninf_{f\in\ball(g,2\gamma_1)\cap F} \phi^\nu((\xi^1, \dots , \xi^\nu),f) \leq \phi^\nu((\xi^1, \dots , \xi^\nu),f_{k^\star}^\star) \leq \phi(f_{k^\star}^\star) + \gamma_2 \leq \phi(g) + \delta - \delta_0.
\]
Thus, the first condition in Proposition \ref{prop:argmin} holds with $\gamma
= 2\gamma_1$ and $\tau = \delta - \delta_0$, and
\[
\exs\big(\nlev_{\leq \delta_0} \phi, \nlev_{\leq \delta_0 + \tau} \phi^\nu((\xi^1, \dots, \xi^\nu),\cdot)\big)
\leq 2\gamma_1.
\]
Equivalently,
\[
\exs\big(\nargmin_{f\in F} \phi(f), \nlev_{\leq \delta} \phi^\nu((\xi^1, \dots , \xi^\nu),\cdot)\big)
\leq 2\gamma_1.
\]

By Theorem \ref{thm:entropy}, there exists $\bar \epsilon>0$ such that $\log
N$ is bounded from above by a term proportional to $\gamma_1^{-n}(\log
\gamma_1^{-1})^{n+1}$ for all $\gamma_1\in (0,\bar\epsilon]$. Thus, there is
a constant $c_1>0$ such that
\begin{equation}\label{eqn:rate}
\frac{\log N - \log(\alpha/2)}{c_0 \gamma_2^2} \leq c_1 \gamma_1^{-n}(\log \gamma_1^{-1})^{n+1} \mbox{ for } \gamma_1 \in (0,\bar\epsilon].
\end{equation}
In view of \eqref{eqn:samplesize}, the right-hand side of \eqref{eqn:rate} provides the rate of increase in sample size that is
needed to guarantee an excess of at most $2\gamma_1$ with confidence level
$1-\alpha$. Inverting the expression, we find that $\gamma_1$ can be
propositional to $\nu^{-1/n} (\log \nu)^{1+1/n}$ as long as $\nu$ is
sufficiently large, which establishes the conclusion.\eop

When $f^\star\in \nargmin_{f\in F} \Ex[\psi(\bfxi^1,f)]$,  the theorem
guarantees that with probability $1-\alpha$
\[
\dist\Big(f^\star, \nlev_{\leq \delta} \Big\{\frac{1}{\nu}\sum_{j=1}^\nu \psi(\bfxi^j, \cdot)\Big\}\Big)
\leq \frac{c (\log \nu)^{1+1/n}}{\nu^{1/n}}
\]
for sufficiently large $\nu$. Hence, the minimizer $f^\star$ of (FIP-U) is
covered by the given level set when appropriately enlarged with a quantity
that vanishes with increasing sample size at nearly the rate $\nu^{-1/n}$.
The confidence region is not given in terms of minimizers of the
approximating problem (FIP-U)$^\nu$, but rather certain level sets.
Membership in such a level set is trivially assessed, does not require
solving the approximating problem, and can be used to rule out the optimality of a candidate $f$. In general, minimizers of (FIP-U)$^\nu$
are not well behaved and depend on the conditioning of (FIP-U) as discussed
in \cite{Royset.18}. Theorem \ref{thm:confidence} bypasses this issue by
considering level sets. Other strengths of Theorem \ref{thm:confidence} are
its mild assumption on the (random) objective function $\psi$ and the wide
range of constraints that is permitted; the family $F$ can be any bounded
closed set in $(\uscfcns(S),\setd)$. The functions $f\mapsto \psi(\xi,f)$ is
only required to be lsc. The assumption about sub-exponential distribution of
$\psi(\bfxi^1,f)$ can be checked pointwise for each $f\in F$. Actually, this
assumption can be relaxed because the proof of Theorem \ref{thm:confidence}
only requires that sample averages are sufficiently low relative to the
actual expectations, but this merely improves $c$ in the theorem and we omit
this refinement.

The practical construction of confidence regions is hampered by the unknown
and hard-to-estimate constants $c$ and $\bar \nu$ in Theorem
\ref{thm:confidence}. In practice, coverage may therefore only be guaranteed
asymptotically. The other unknown parameter $\delta$ is easy to estimate
conservatively because for any $f\in F$, the sample average
$\frac{1}{\nu}\sum_{j=1}^\nu \psi(\bfxi^j,f)$, using a different sample,
furnishes an estimator of $\Ex[\psi(\bfxi^1,f)]$, which in turn is an upper
bound on $\ninf_{f\in F} \Ex[\psi(\bfxi^1,f)]$. An effort to select a low
$\delta$ would obviously result in a smaller level set, but typically also
large $c$ and $\bar \nu$.

The effect of $n$ on the rate of convergence is profound and in line with the
growth of the covering numbers as $n$ increase. It highlights, for example,
the fundamental challenge associated with high-dimensional nonparametric
estimation already well documented (see
\cite{BalabdaouiWellner.10,KimSamworth.16}). On the positive note, if $n=1$,
which already captures many interesting applications \cite{RoysetWets.18a},
then the convergence rate is nearly $\nu^{-1}$ and therefore {\it faster}
than the canonical $\nu^{-1/2}$ rate. If $F$ is restricted to some
finite-dimensional subset of $\uscfcns(S)$, then the covering numbers from
Theorem \ref{thm:entropy} can be replaced by much improved ones, typically of
order $O(\epsilon^{-1})$ so that their logarithm is of order $O(\log
\epsilon^{-1})$ and the rate improves from essentially $\nu^{-1/n}$ to
$e^{-\nu}$ in Theorem \ref{thm:confidence}.

We illustrate the application of Theorem \ref{thm:confidence} on stochastic
optimization problems arising in nonparametric statistics. \\

\state Example 1: Maximum Likelihood Estimation of Probability Densities.
Suppose that we would like to estimate an unknown probability density
function $f^0\in\uscfcns(S)$. Since we permit densities to have value zero on
a subset of $S$, there is no requirement that the support of $f^0$ is known;
$S$ just needs to contain the support. Given a sample $\bfxi^1, \dots,
\bfxi^\nu$, which in this case takes values in $S$, i.e., $\Xi = S$, a
maximum likelihood estimator of $f^0$ over a class $F\subset \uscfcns(S)$ is
any minimizer of
\[
\nmin_{f\in F} -\frac{1}{\nu}\sum_{j=1}^\nu \log f(\bfxi^j)
\]
and, in the notation above, $\psi(\xi,f) = -\log f(\xi)$. The function
$(\xi,f)\mapsto -\log f(\xi)$ is a random lsc function on the probability
space $(S,\cB,P)$, where $P$ is the probability distribution of $f^0$ and
$\cB$ contains the Borel sets of $(S,\|\cdot-\cdot\|_\infty)$ supplemented
with the necessary probability-zero sets to make the probability space
complete. This fact is easily realized because the function is actually lsc
jointly in its arguments; see \cite{RoysetWets.18} for details.

In this case, $\psi(\bfxi^1,f)$ being sub-exponential amounts to having $F$
consist of sub-exponential densities. The requirement about inf-integrability
is extensively discussed in \cite{RoysetWets.18}. For example, suppose $F$ is
a nonempty closed subset of
\[
\Bigg\{f\in \uscfcns(S)~\Bigg| ~\int f(x) dx = 1, ~\int x f(x) dx \in C, ~u(x) \leq f(x) \leq v(x), ~\forall x\in S\Bigg\},
\]
where $C\subset \reals^n$ is closed and $u,v:S\to (0,\infty)$, with $v\in
\uscfcns(S)$. Moreover, suppose that the actual density $f^0\in F$ and for
some $\gamma_1\geq 0,\gamma_2>0$, and, $\zeta_1,\zeta_2\in \reals$,
\[
u(x) \geq e^{-\gamma_1 \|x\|_\infty + \zeta_1} \mbox{ and } v(x) \leq e^{-\gamma_2 \|x\|_\infty + \zeta_2}.
\]
All the assumptions of Theorem \ref{thm:confidence} are then satisfied. The
requirement that $F$ is bounded is automatically satisfied because $f\geq 0$
for all $f\in F$.

In this example the maximum likelihood estimator finds the best estimate that
satisfies the given pointwise bounds and moment restriction. There is {\it
no} requirement that the actual density or its estimate should be smooth or
even continuous. Of course, a large variety of other constraints can be
brought in too; see \cite{RoysetWets.18} for some possibilities.\\

We recall that a subset $F_0 \subset \uscfcns(S)$ is {\it equi-usc}
\cite[Sect. 7.B]{VaAn} if there exists $\delta:S\times (0,\infty) \times
(0,\infty) \to (0,\infty)$ such that for any $\epsilon,\rho>0$, $\bar x\in
S$, and $f\in F_0$,
\[
\nsup_{x\in \ball_\infty(\bar x,\delta(\bar x,\epsilon,\rho))} f(\bar x) \leq \max\{f(x) + \epsilon, -\rho\}.
\]
If $F_0$ is a singleton, then the condition reduces to that of usc. If $F_0$
contains only Lipschitz continuous functions, or only piecewise Lipschitz
continuous functions, or only finite-valued concave functions on $\reals^n$,
to mention some examples, then $F_0$ is equi-usc.\\

\state Example 2: Least-Squares Regression. Suppose that we are given the
random design model
\[
\bfy^j = f^0(\bfx^j) + \bfz^j, ~~j=1, 2, \dots, \nu
\]
where $\bfx^1, \bfx^2, \dots, \bfx^\nu$ are independent and identically
distributed $n$-dimensional random vectors that take values in a closed set
$S\subset\reals^n$,  $\bfz^1, \bfz^2, \dots, \bfz^\nu$ are zero-mean random
variables that are also independent of $\bfx^1, \bfx^2, \dots, \bfx^\nu$, and
$f^0:S\to\reals$ is an unknown function to be estimated based on observations
of $\bfxi^1 = (\bfx^1,\bfy^1)$. In this case, $\Xi = S \times \reals$, again
we adopt a sigma-algebra that contains the Borel sets on $\Xi$ and that
results in a complete probability space under the distribution of
$(\bfx^1,\bfy^1)$. The least-squares estimator of $f^0$ over the class
$F\subset\uscfcns(S)$ is then any minimizer of
\[
\nmin_{f\in F} \frac{1}{\nu}\sum_{j=1}^\nu \big(\bfy^j - f(\bfx^j)\big)^2.
\]
Resulting estimates furnish approximations of $f^0$ that in an engineering
design context can be maximized to find an optimal design without any
(additional) costly simulation of system performance. The only simulations
required are those needed to generate a data set $\{(x^j,y^j), j=1, ...,
\nu\}$.

In this case, $\psi((x,y),f) = (y - f(x))^2$. Since $(x,f)\mapsto f(x)$ is
usc and thus measurable, we also have that $\psi$ is measurable.
Consequently, $\psi$ is a random lsc function provided that $F$ is equi-usc,
an assumption that provides the necessary pointwise convergence (cf.
\cite[Thm. 7.10]{VaAn}). Its nonnegativity ensures that $\psi$ is also
locally inf-integrable.

A confidence region for $f^0\in F$ emerges from Theorem \ref{thm:confidence}
when $(\bfy^1 - f(\bfx^1))^2$ is sub-exponential for all $f\in F$. For
example, this will be the case when $\bfz^1$ and every component of $\bfx^1$
are sub-Gaussian, and for some $\gamma,\zeta\in\reals$,
\[
f\in F \Longrightarrow |f(x)| \leq \gamma\|x\|_\infty + \zeta, ~\forall x\in S.
\]
Since $f^0$ must be a minimizer of $\nmin_{f\in F} \Ex[(\bfy^1 -
f(\bfx^1))^2]$, provided that $f^0\in F$, Theorem \ref{thm:confidence}
guarantees that $f^0$ is covered by the stipulated level set when
appropriately enlarged.

\subsection{Rates of Convergence under H\"{o}lder Condition}

Theorem \ref{thm:confidence} does not rule out the possibility that the limit
of the given level sets strictly contains $\nargmin_{f\in F}
\Ex[\psi(\bfxi^1,f)]$. In fact, this cannot be ruled out unless additional
assumptions are brought in; \cite{Royset.18} contains a discussion. Still, a
H\"{o}lder condition enables us to ``reverse'' Theorem \ref{thm:confidence}
and quantify the rate of convergence of the excess of minimizers of
(FIP-U)$^\nu$ over those of (FIP-U). Since it is relatively straightforward,
we also address approximating constraints. Although the approximating constraints can be rather general, the rate of convergence in the following theorem depends on the rate with which the approximating feasible set approaches the actual one. Thus, it is not immediately clear how the piecewise affine functions discussed in Section 3, which have unknown rate of convergence, can be used for constructing these approximations.

As usual, we let $\alpha^\nu =
o(r^\nu)$ imply that for every $\delta>0$ there exists $\bar \nu$ such that
$\alpha^\nu \leq \delta r^\nu$ for all $\nu\geq \bar \nu$.

\begin{theorem}{\rm (rate of convergence).}\label{thm:rate} For a complete probability space $(\Xi,\cA,\mathbb{P})$ and closed and
bounded sets $F^\nu,F^0\subset F\subset\uscfcns(S)$, suppose that
$\psi:\Xi\times F\to (-\infty,\infty]$ is a random lsc function for which
there exist $p\in (0,\infty)$ and integrable random variable $\kappa:\Xi\to
[0,\infty)$ such that
\[
|\psi(\xi,f) -  \psi(\xi,g)| \leq \kappa(\xi) [ \setd(f,g) ]^p \mbox{ for all } f,g\in F \mbox{ and } \xi\in \Xi.
\]
Suppose also that $\bfxi^1, \bfxi^2, \dots$ are independent random elements,
each with values in $\Xi$ and distributed according to $\mathbb{P}$, and
$\psi(\bfxi^1,f)$ is sub-exponential for all $f\in F$. Let
\[
r^\nu = \nu^{\frac{-1}{2+n/p}}(\log \nu)^{\frac{1+n}{2+n/p}}.
\]
If $\mathbb{H}(F^\nu,F^0) = o(\min\{r^\nu,(r^\nu)^{1/p}\})$ and $\alpha\in
(0,1)$, then there exist $c\in [0,\infty)$ and $\bar \nu\in\nats$ such that
for $\nu\geq \bar \nu$ and $\epsilon^\nu\geq 0$,
\[
P^\nu\bigg[\exs\Big(\epsilon^\nu\mbox{-}\nargmin_{f\in F^\nu} \frac{1}{\nu}\sum_{j=1}^\nu \psi(\bfxi^j, f), ~~(\epsilon^\nu+c r^\nu)\mbox{-}\nargmin_{f\in F^0} \Ex[\psi(\bfxi^1,f)]\Big)
\leq \mathbb{H}(F^\nu,F^0)\bigg] \geq 1 - \alpha.
\]
\end{theorem}
\state Proof. Let $\zeta
>0$. Since $\kappa(\bfxi^1)$ is integrable, there exists $\bar \nu_0\in\nats$ such that
$P^\nu(|\nu^{-1}\sum_{j=1}^\nu \kappa(\bfxi^j) - \Ex[\kappa(\bfxi^1)]| \geq
\zeta ) \leq \alpha/2$ for all $\nu\geq \bar \nu_0$. Let $\gamma_1>0$. As in
Theorem \ref{thm:entropy}, there is a finite number $N = N(F,\gamma_1/2)$ of
closed balls in $(\uscfcns(S),\setd)$ with radius $\gamma_1/2$ and center
$f_k'$ covering $F$. To make sure that the balls are centered at points in
$F$, we can always select some other centers $f_k\in F$ and balls with radius
$\gamma_1$ and still cover $F$.

Let $\phi$ and $\phi^\nu$ be as defined in the proof of Theorem
\ref{thm:confidence}. We note that $\psi$ is locally inf-integrable due to
the H\"{o}lder condition and the pointwise sub-exponential property. Since
$\psi(\bfxi^1,f_k)$ is sub-exponential, there exists by Bernstein's
inequality $\bar\gamma_2>0$ and $c_0>0$ such that for $\gamma_2 \in [0,
\bar\gamma_2]$,
\[
P^\nu\big(|\phi^\nu((\bfxi^1, \dots , \bfxi^\nu), f_k) - \phi(f_k)| \geq \gamma_2\big) \leq 2e^{-\nu c_0\gamma_2^2} \mbox{ for all } k=1, \dots , N.
\]
Consequently, as long as $\nu\geq \bar \nu_0$ and $2Ne^{-\nu c_0\gamma_2^2}
\leq \alpha/2$, or, equivalently,
\[
\nu\geq \max\Big\{\bar\nu_0, \frac{\log N - \log(\alpha/4)}{c_0 \gamma_2^2}\Big\}
\]
we have that
\[
P^\nu\bigg(\max_{k=1, \dots ,N} \Big|\phi^\nu((\bfxi^1,\dots , \bfxi^\nu), f_k) - \phi(f_k)\Big| \geq \gamma_2 ~\mbox{ or } ~\Big|\frac{1}{\nu}\sum_{j=1}^\nu \kappa(\bfxi^j) - \Ex[\kappa(\bfxi^1)]\Big| \geq \zeta\bigg) \leq \alpha.
\]
Suppose that we have an event $(\xi^1, \dots , \xi^\nu)\in \Xi^\nu$ where
\[
\max_{k=1, \dots ,N} \Big|\phi^\nu((\xi^1,\dots , \xi^\nu), f_k) - \phi(f_k)\Big| < \gamma_2 ~\mbox{ and }~ \Big|\frac{1}{\nu}\sum_{j=1}^\nu \kappa(\xi^j) - \Ex[\kappa(\bfxi^1)]\Big| < \zeta.
\]

Next, we apply Proposition \ref{prop:argmin} for the lsc functions $\bar
\phi:F^0\to \reals$ given by $\bar \phi(f) = \phi(f)$ and $\bar
\phi^\nu:F^\nu\to \reals$ given by $\bar\phi^\nu(f) = \phi^\nu((\xi^1, \dots,
\xi^\nu),f)$. In view of the H\"{o}lder assumption on $\psi$, this implies
that $\bar\phi^\nu$ is finite when defined. Moreover, for all $f,g\in F$,
\[
\big|\phi(f) - \phi(g)\big| \leq \Ex[\kappa(\bfxi^1)] [\setd(f,g)]^p.
\]
Let $\delta^\nu = \mathbb{H}(F^\nu,F^0)$. Suppose that $f\in F^\nu$. Then,
there is $f'\in F^0$ and $k^\star \in \{1, \dots, N\}$ such that
$\setd(f,f')\leq \delta^\nu$ and $\setd(f',f_{k^\star}) \leq \gamma_1$. Thus,
\begin{align*}
&\ninf_{g\in \ball(f,\delta^\nu)\cap F} \bar\phi(g) \leq  \bar\phi(f') \leq  \bar\phi(f_{k^\star}) + \Ex[\kappa(\bfxi^1)] \gamma_1^p\\
&<     \bar\phi^\nu(f_{k^\star}) + \gamma_2 + \Ex[\kappa(\bfxi^1)] \gamma_1^p\\
&\leq  \bar\phi^\nu(f) + \gamma_2 + \Ex[\kappa(\bfxi^1)] \gamma_1^p +  \frac{1}{\nu}\sum_{j=1}^\nu \kappa(\xi^j)(\delta^\nu + \gamma_1)^p\\
&\leq  \bar\phi^\nu(f) + \gamma_2 + \Ex[\kappa(\bfxi^1)](\gamma_1^p + (\delta^\nu + \gamma_1)^p) + \zeta(\delta^\nu + \gamma_1)^p.
\end{align*}
Similarly, suppose that $f\in F^0$. Then, there is $f'\in F^\nu$ and $k^\star
\in \{1, \dots, N\}$ such that $\setd(f,f')\leq \delta^\nu$ and
$\setd(f',f_{k^\star}) \leq \gamma_1$. Consequently,
\begin{align*}
&\ninf_{g\in \ball(f,\delta^\nu)\cap F} \bar\phi^\nu(g)\leq  \bar\phi^\nu(f') \leq  \bar\phi^\nu(f_{k^\star}) + \frac{1}{\nu}\sum_{j=1}^\nu \kappa(\xi^j)\gamma_1^p\\
&< \bar\phi(f_{k^\star}) + \gamma_2  + \frac{1}{\nu}\sum_{j=1}^\nu \kappa(\xi^j)\gamma_1^p\\
&\leq  \bar\phi(f) + \Ex[\kappa(\bfxi^1)](\delta^\nu + \gamma_1)^p + \gamma_2  + \frac{1}{\nu}\sum_{j=1}^\nu \kappa(\xi^j)\gamma_1^p\\
&\leq  \bar\phi(f) + \Ex[\kappa(\bfxi^1)][(\delta^\nu + \gamma_1)^p + \gamma_1^p] + \gamma_2  + \zeta\gamma_1^p.
\end{align*}
Thus, we have shown that the conditions of Proposition \ref{prop:argmin} hold
for the functions $\bar\phi$ and $\bar\phi^\nu$ with
\[
\delta^\nu \mbox{ and } \tau_0 = \gamma_2 + \Ex[\kappa(\bfxi^1)](\gamma_1^p + (\delta^\nu + \gamma_1)^p) + \zeta(\delta^\nu + \gamma_1)^p
\]
as the two error parameters ($\gamma$ and $\tau$) and we therefore have that
\[
\exs\big(\epsilon^\nu\mbox{-}\nargmin_{f\in F^\nu}
\phi^\nu((\xi^1, \dots , \xi^\nu), f), ~(\epsilon^\nu+2\tau_0)\mbox{-}\nargmin_{f\in F^0} \phi(f) \big) \leq
\delta^\nu.
\]
By Theorem \ref{thm:entropy}, $\log N$ is bounded from by a term proportional
to $\gamma_1^{-n}(\log \gamma_1^{-1})^{n+1}$ for sufficiently small
$\gamma_1$. Thus, there exist constants $c_1,c_2>0$ such that
\[
\frac{\log N - \log(\alpha/4)}{c_0 \gamma_2^2} \leq c_1 \gamma_1^{-n}(\log \gamma_1^{-1})^{n+1}\gamma_2^{-2} + c_2 \gamma_2^{-2},
\]
which gives the rate of growth in $\nu$ as $\gamma_1$ and $\gamma_2$ vanish.
For $\tau>0$, the error $\tau_0$ can be kept below $\tau$ if $\gamma_1$ is
proportional to $\tau^{1/p}$, $\gamma_2$ is proportional to $\tau$,
$\delta^\nu$ is proportional to $\min\{\tau, \tau^{1/p}\}$, and the
(positive) proportionality constants are selected sufficiently close to zero.
In view of these choices about $\gamma_1$ and $\gamma_2$, there is a constant
$c_3>0$ such that
\[
c_1 \gamma_1^{-n}(\log \gamma_1^{-1})^{n+1}\gamma_2^{-2} + c_2 \gamma_2^{-2} \leq c_3 \tau^{-2-n/p}(\log \tau^{-1/p})^{n+1}.
\]
With $\nu$ above $\bar \nu_0$ as well as the previous right-hand side, or
equivalently for some $c_4>0$,
\[
\tau \geq c_4 \nu^{\frac{-1}{2+n/p}} (\log \nu)^{\frac{1+n}{2+n/p}},
\]
we ensure the required confidence level and the conclusion follows.\eop

A corollary of the theorem for the case with $\epsilon^\nu =0$ and $F^\nu =
F^0 = F$ is that
\[
\nargmin_{f\in F} \frac{1}{\nu} \sum_{j=1}^\nu \psi(\bfxi^j,f) ~ \subset ~cr^\nu\mbox{-}\nargmin_{f\in F} \Ex[\psi(\bfxi^1,f)]
\]
with at least probability $1-\alpha$. Thus, minimizers of (FIP-U)$^\nu$
converge at the rate $r^\nu$ to a minimizer of (FIP-U). The rate depends on
the H\"{o}lder coefficient $p$ as well as the dimension $n$ of the space of
function under considerations.

We illustrate the assumptions of the theorem for two stochastic optimization
problems arising in nonparametric statistics, but start with an intermediate
result.
\begin{proposition}\label{prop:ptwise}
  For Lipschitz continuous functions $f,g\in \uscfcns(S)$
   with common modulus $\kappa\in
  [0,\infty)$,
  \[
|f(x) - g(x)| \leq (1+\kappa)e^{\rho(x)} \setd(f,g) \mbox{ for all } x\in S,
  \]
where $\rho(x) = \max\{\|x\|_\infty, |f(x)|, |g(x)|\}$.
\end{proposition}
\state Proof. Let $x\in S$. The first result is trivial if $\rho(x) =
\infty$. Suppose that $\rho(x)<\infty$. From Lemma \ref{lemma:distestim},
$\setd(f,g) \geq e^{-\rho(x)} \hatsetd_{\rho(x)}(f,g)$. Set $\tau \in
(\hatsetd_{\rho(x)}(f,g),\infty)$. Again, by Lemma \ref{lemma:distestim},
there exists $y\in \ball_\infty(x,\tau)$ such that $f(y) \geq g(x) - \tau$.
Thus, $g(x) - f(x) = g(x) - f(y) + f(y) - f(x) \leq \tau + \kappa\tau$. A
similar argument establishes that $f(x) - g(x) \leq \tau + \kappa\tau$.
Hence, by letting $\tau$ tends to its lower limit, we obtain that $|f(x) -
g(x)| \leq (1 + \kappa)\hatsetd_{\rho(x)}(f,g)$ and the conclusion follows.
\eop\\

\state Example 3: Least-Squares Regression. We return to the setting of
Example 2, but now let $F$ be a family that contains only Lipschitz
continuous functions with common modulus $\kappa_0\geq 0$. Suppose also that
$\bfz^1$ and every component of $\bfx^1$ are sub-Gaussian, the unknown
function $f^0\in F$, and there exists $\beta<\infty$ such that $f(0)\leq
\beta$ for all $f\in F$. Then, $F$ is equi-usc and there are
$\gamma,\zeta\in\reals$ such that $|f(x)| \leq \gamma\|x\|_\infty+\zeta$ for
all $f\in F$. Proposition \ref{prop:ptwise} then ensures that the H\"{o}lder
condition in Theorem \ref{thm:rate} holds with $p=2$ and $\kappa((x,y)) =
(1+\kappa_0)^2 \exp(2\max\{\|x\|_\infty, \gamma\|x\|_\infty +\zeta\})$, which
is integrable in view of the sub-Gaussianity of $\bfx^1$.

 Using the
bound on $|f(x)|$, we also have that $(\bfy^1-f(\bfx^1))^2 =
(f^0(\bfx^1)-f(\bfx^1) + \bfz^1)^2$ is sub-exponential. The assumptions of
Theorem \ref{thm:rate} therefore hold,
\[
r^\nu = \nu^{\frac{-2}{4+n}}(\log \nu)^{\frac{1+n}{2+n/2}},
\]
and, for closed $F^\nu,F^0\subset F$, there exist $c\in [0,\infty)$ and $\bar
\nu\in\nats$ such that
\[
P^\nu\bigg[\exs\Big(\nargmin_{f\in F^\nu} \frac{1}{\nu}\sum_{j=1}^\nu (\bfy^j - f(\bfx^j))^2, ~~cr^\nu\mbox{-}\nargmin_{f\in F^0} \Ex[(\bfy^1 - f(\bfx^1))^2]\Big) \leq \mathbb{H}(F^\nu,F^0)\bigg] \geq 1 - \alpha
\]
provided that $\nu\geq \bar \nu$ and $\mathbb{H}(F^\nu,F^0) = o(r^\nu)$.
Thus, when $\mathbb{H}(F^\nu,F^0) =0$ and $\hat f^\nu\in \nargmin_{f\in
F^\nu} \frac{1}{\nu}\sum_{j=1}^\nu$ $(\bfy^j - f(\bfx^j))^2$ is measurable,
\[
P^\nu\Big(\Ex\big[(\hat f^\nu(\bfx^1) - f^0(\bfx^1))^2\big] \leq cr^\nu\Big) \geq 1-\alpha.
\]
The rates developed here apply in rather general settings and remain in
effect even if $f^0\not\in F^0$. More specific settings give improved results
as in the case of regression with fixed design and Lipschitz continuous
functions defined on compact convex subset \cite[p.
333]{vanderVaartWellner.96} and in the univariate case
\cite{GuntuboyinaSen.15}.\\

\state Example 4: Least-Squares Probability Density Estimation. We return to
the setting of Example 1, but now consider the least-squares estimator of
$f^0$, which is any minimizer of
\[
\nmin_{f\in F^\nu} -\frac{2}{\nu}\sum_{j=1}^\nu f(\bfxi^j) + \int [f(x)]^2 dx.
\]
This estimator is motivated by the fact that the unknown function
\[
f^0 \in \nargmin_{f\in F} \int \big[f(x) - f^0(x)\big]^2 dx = \nargmin_{f\in F} -2\Ex[f(\bfxi^1)] + \int [f(x)]^2 dx
\]
whenever $f^0\in F$. To make the case rather concrete, let $\kappa\in
[0,\infty)$ and for some bounded function $h:S\to [0,\infty)$, with $\int
h(x) dx < \infty$,
\[
F = \bigg\{f\in \uscfcns(S)~\bigg|~\int f(x) dx = 1, ~0 \leq f(x) \leq h(x), ~|f(x) - f(y)|\leq \kappa\|x-y\|_\infty, \forall x,y\in S\bigg\},
\]
which can be shown to be closed and bounded; see arguments in
\cite{RoysetWets.18}. In this case, $(\xi,x) \mapsto \psi(\xi,f) = -2f(\xi) +
\int [f(x)]^2 dx$ is a random lsc function as can be seen by invoking Fatou's
Lemma and pointwise convergence; again see \cite{RoysetWets.18}.  Then,
$\psi(\bfxi^1,f)$ is sub-exponential for all $f\in F$ as it is in fact
bounded.

It remains to check the H\"{o}lder condition in Theorem \ref{thm:rate}.
Suppose that $\Ex[\exp(\|\bfxi^1\|_\infty)]<\infty$. In view of Proposition
\ref{prop:ptwise}, if the integral term in $\psi$ had not been present, then
the condition holds with $p=1$; Lipschitz continuity and the fact that
$\Ex[\exp(\|\bfxi^1\|_\infty)]<\infty$ ensures integrability of the
H\"{o}lder modulus. If $S$ were compact, then $\psi$ would still satisfy the
condition with $p=1$. For a noncompact $S$, the argument needs to be slightly
modified by first ``ignoring'' the integral term and second reintroduce it in
a slightly generalized version of Theorem \ref{thm:rate}. We omit the
details.

In summary, for the given $F$ and under the assumption that
$\Ex[\exp(\|\bfxi^1\|_\infty)]<\infty$, we can show by invoking Theorem
\ref{thm:rate} (or the mentioned extensions) that for any $\alpha\in (0,1)$
there exist $c\in [0,\infty)$ and $\bar \nu\in\nats$ such that for every
$\nu\geq \bar \nu$ and $\epsilon^\nu\geq 0$
\begin{align*}
& P^\nu\bigg[\epsilon^\nu\mbox{-}\nargmin_{f\in F}
 -\frac{2}{\nu}\sum_{j=1}^\nu f(\bfxi^j) + \int [f(x)]^2 dx \subset\\
&~~~~~~~~~~~~~~(\epsilon^\nu+cr^\nu)\mbox{-}\nargmin_{f\in F} -2\Ex[f(\bfxi^1)] + \int [f(x)]^2
dx\bigg] \geq 1 - \alpha \mbox{ with } r^\nu = \nu^{\frac{-1}{2+n}}(\log \nu)^{\frac{1+n}{2+n}}.
\end{align*}
A sharper result is available in the univariate case over the class of
nonincreasing convex functions \cite{GroeneboomJongbloedWellner.01}.\\

\state Acknowledgements. This work in supported in parts by DARPA under
grants HR0011-14-1-0060 and HR0011-8-34187, and Office of Naval Research
(Science of Autonomy Program) under grant N00014- 17-1-2372.

\section*{Appendix}

\noindent {\bf Proof of Theorem \ref{thm:mainlow}.} Let $\rho>0$ and $F =
\{f\in\uscfcns(\reals^d)~|~f(x)\geq -\rho \mbox{ for at least one } x\in
[0,\rho]^n\}$. We show that $F$ cannot be covered with a lower number of
balls than stipulated. Clearly, $\dist_\infty(0,\hypo f) \leq \rho$ for all
$f\in F$. Thus, in view of (\ref{eqn:distbds}), $\setd(0,f) \leq \rho + 1$
for all $f\in F$, where $0$ is the zero function on $\reals^n$, and $F$ is
therefore bounded.

Next, let $\epsilon \in (0,\rho e^{-\rho}/6]$. We discretize $[0,\rho]^n$ by
defining $x_i^k = k \rho/\nu_\epsilon$, $k = 1, ..., \nu_\epsilon-1$ and
$i=1, ..., n$, where
\begin{equation*}
  \nu_\epsilon = \left\lfloor \frac{\rho e^{-\rho}}{3\epsilon}\right\rfloor \geq 2,
\end{equation*}
with $\lfloor a \rfloor$ being the largest integer not exceeding $a$. The
discretization of $[0,\rho]^n$ then contains the points $(x_1^{k_1},
x_2^{k_2}, ..., x_n^{k_n})$, with $k_i \in \{1, 2, ..., \nu_\epsilon-1\}$ and
$i=1, ..., n$. Clearly, the distance between any two such points in the
sup-norm is at least $\rho/\nu_\epsilon\geq 3\epsilon e^\rho$. We carry out a
similar discretization of $[-\rho,0]$ and define $y^l = l \rho/
\nu_\epsilon$, $l=1, ..., \nu_\epsilon$. The functions that are finite on the
discretization points of $[0,\rho]^n$, with values at each such point equal
to $y^l$ for some $l$, and have value minus infinity elsewhere are given by
$F_{\epsilon}$, i.e.,
\begin{align*}
  F_{\epsilon} =& \{f\in \uscfcns(\reals^n)~|~ \mbox{for each } x=(x_1^{k_1}, ..., x_n^{k_n}), \mbox{ with } k_i \in \{1, 2, ..., \nu_\epsilon-1\}, f(x) = y^l\\
  &   \mbox{ for some } l=1, ..., \nu_\epsilon; f(x) = -\infty \mbox{ otherwise} \}.
\end{align*}
Certainly, $F_{\epsilon} \subset F$. We next define
\begin{equation*}
  G_\epsilon(f) = \{g\in \uscfcns(\reals^n)~|~ \hatsetd_\rho(f,g) \leq \epsilon e^\rho\}, ~\mbox{ for } f\in\uscfcns(\reals^n).
\end{equation*}
We establish that $G_\epsilon(f) \cap G_\epsilon(f') = \emptyset$ for
$f,f'\in F_{\epsilon}, f\neq f'$. Suppose for the sake of a contradiction
that there is a $g$ with $g\in G_\epsilon(f)$ and $g\in G_\epsilon(f')$ for
$f,f'\in F_\epsilon$, $f\neq f'$. Then, $\hatsetd_\rho(f,g) \leq \epsilon
e^\rho$ and $\hatsetd_\rho(f',g) \leq \epsilon e^\rho$. However, since $f\neq
f'$, there exists a point $x\in [0,\rho]^n$ with $|f(x) - f'(x)| \geq 3
\epsilon e^\rho$. Without loss of generality, suppose that $f(x) \geq f'(x) +
3\epsilon e^\rho$. Since $f(z), f'(z) = -\infty$ for all $z\neq x$ with
$\|z-x\|_\infty < 3\epsilon e^\rho$, we have that $\hatsetd_\rho(f,g) \leq
\epsilon e^\rho$ implies that $g(z) \geq f(x) - \epsilon e^\rho$ for some
$z\in \ball(x,\epsilon e^\rho)$. Moreover, $\hatsetd_\rho(f',g) \leq \epsilon
e^\rho$ implies that $g(z) \leq f'(x) + \epsilon e^\rho \leq f(x) - 3\epsilon
e^\rho + \epsilon e^\rho = f(x) - 2\epsilon e^\rho$ for all $z\in
\ball(x,\epsilon e^\rho)$. Since this is not possible for $g$, we have
reached a contradiction. Thus, $G_\epsilon(f) \cap G_\epsilon(f') =
\emptyset$ for $f,f'\in F_{\epsilon}, f\neq f'$.

By Lemma \ref{lemma:distestim}, for any $f\in \uscfcns(\reals^n)$,
\begin{equation*}
  \setd(f,g) \geq e^{-\rho} \hatsetd_\rho(f,g) > e^{-\rho} \epsilon e^\rho = \epsilon \mbox{ for all } g\not\in G_\epsilon(f).
\end{equation*}
Hence, for $f\in F_{\epsilon}$, an $\setd$-ball with radius $\epsilon$ that
contains $f$ needs to be centered at some $g\in G_\epsilon(f)$. Since the
sets $G_\epsilon(f)$, $f\in F_{\epsilon}$, are nonoverlapping, a cover of
$F_{\epsilon}$ by $\setd$-balls with radius $\epsilon$ must involve a number
of balls that is at least as great as the number of functions in
$F_{\epsilon}$, which is $\nu_\epsilon^{m_\epsilon}$, where $m_\epsilon =
(\nu_\epsilon-1)^n$. Thus,
\begin{equation}\label{eqn:coveringLB}
  \log N(F,\epsilon) \geq \nu_\epsilon^n \log \nu_\epsilon \geq \left(\frac{\rho e^{-\rho}}{3\epsilon}-2\right)^n \log \left(\frac{\rho e^{-\rho}}{3\epsilon}-1\right).
\end{equation}
Let $c_1 = |\log (\rho e^{-\rho}/4)|$ and $\bar\epsilon = \min\{\rho
e^{-\rho}/12, e^{-2c_1}\}$. Continuing from \eqref{eqn:coveringLB}, we then
find that
\begin{equation*}
  \log N(F,\epsilon) \geq \left(\frac{\rho e^{-\rho}}{6}\right)^n \left[1+ \frac{\log(\rho e^{-\rho}/4)}{\log \epsilon^{-1}} \right] \frac{1}{\epsilon^n}\log \frac{1}{\epsilon}.
\end{equation*}
Since $\log \epsilon^{-1} \geq 2|\log (\rho e^{-\rho}/4)|$ for $\epsilon \in
(0, \bar \epsilon]$, we have that
\begin{equation*}
  \log N(F,\epsilon) \geq \left(\frac{\rho e^{-\rho}}{6}\right)^n \frac{1}{2}\frac{1}{\epsilon^n}\log \frac{1}{\epsilon}~ \mbox{ for } \epsilon \in (0, \bar \epsilon],
\end{equation*}
and the conclusion is reached.\eop

\bibliographystyle{plain}
\bibliography{refs}

\end{document}